\documentclass[12pt]{article}
\usepackage[final]{epsfig}
\usepackage{graphics}
\usepackage{amsmath}
\usepackage{amsfonts}
\usepackage{latexsym}
\usepackage{amssymb}
\usepackage{graphicx}
\usepackage{epstopdf}
\newtheorem{lemma}{Lemma}[section]
\newtheorem{proposition}[lemma]{Proposition}
\newtheorem{remark}[lemma]{Remark}
\newtheorem{example}[lemma]{Example}
\newtheorem{theorem}[lemma]{Theorem}

\newtheorem{corollary}[lemma]{Corollary}
\newtheorem{problem}[lemma]{Problem}

\begin{document}
\newcommand{\eps}{{\varepsilon}}
\newcommand{\proofend}{$\Box$\bigskip}
\newcommand{\C}{{\mathbf C}}
\newcommand{\Q}{{\mathbf Q}}
\newcommand{\R}{{\mathbf R}}
\newcommand{\Z}{{\mathbf Z}}
\newcommand{\RP}{{\mathbf {RP}}}

\def\proof{\paragraph{Proof.}}


\newcommand{\marginnote}[1]
{
}

\newcounter{bk}
\newcommand{\bk}[1]
{\stepcounter{bk}$^{\bf BK\thebk}$%
\footnotetext{\hspace{-3.7mm}$^{\blacksquare\!\blacksquare}$
{\bf BK\thebk:~}#1}}

\newcounter{st}
\newcommand{\st}[1]
{\stepcounter{st}$^{\bf ST\thest}$%
\footnotetext{\hspace{-3.7mm}$^{\blacksquare\!\blacksquare}$
{\bf ST\thest:~}#1}}


\title {Pseudo-Riemannian geodesics and billiards}
\author{Boris Khesin\thanks{
Department of Mathematics,
University of Toronto, Toronto, ON M5S 2E4, Canada;
e-mail: \tt{khesin@math.toronto.edu}
}
\, and Serge Tabachnikov\thanks{
Department of Mathematics,
Pennsylvania State University, University Park, PA 16802, USA;
e-mail: \tt{tabachni@math.psu.edu}
}
\\
}
\date{March 22, 2007}
\maketitle
\begin{abstract}
In pseudo-Riemannian geometry the spaces of space-like and  time-like geodesics 
on a pseudo-Riemannian manifold have natural symplectic structures
(just like in the Riemannian case),  while
the space of light-like geodesics has a natural contact structure. 
Furthermore, the space of all geodesics has a structure of a Jacobi manifold.
We describe the geometry of these structures and their generalizations, 
as well as introduce and study pseudo-Euclidean billiards. 
We prove  pseudo-Euclidean analogs of the 
Jacobi-Chasles theorems and show the integrability of 
the billiard in the ellipsoid and the geodesic flow on the ellipsoid
in a pseudo-Euclidean space. 
\end{abstract}

\tableofcontents

\section{Introduction} \label{intro}

The space of oriented lines in the Euclidean $n$-space
has a natural symplectic structure. So does the space 
of geodesics on a Riemannian manifold, at least locally. The structures
on the space of geodesics on a {\it pseudo-Riemannian} manifold
are more subtle. It turns out that the spaces of space-like lines
and  time-like lines in a pseudo-Euclidean space 
have natural symplectic structures (and so do the corresponding
spaces of the geodesics on a pseudo-Riemannian manifold),  while
the space of light-like (or, null) lines or geodesics 
has a natural contact structure. Moreover, the corresponding 
symplectic structures on the manifolds of 
space- and time-like geodesics blow up as one
approaches the border between them, the space of the  null geodesics.
On the other hand, the space of all (space--like, time-like, and null)
geodesics together locally has a structure of a smooth Jacobi manifold.
Below we describe these structures in the pseudo-Riemannian setting, 
emphasizing the differences from the Riemannian case (see Section \ref{str}).

Many other familiar facts in Euclidean/Riemannian geometry have their 
analogs in the pseudo-Riemannian setting, but often with 
an unexpected twist. For example, assign the oriented normal line to each point of a cooriented hypersurface in pseudo-Euclidean space; this gives a smooth map from the hypersurface to the space of oriented lines whose image is Lagrangian in the space of of space-like  and  time-like lines, and Legendrian in the space of light-like  lines, see Section \ref{hyps}.  Another example:   a convex hypersurface
in Euclidean space $\R^n$ has at least $n$ diameters. It turns out 
that a convex hypersurface in pseudo-Euclidean space $V^{k+l}$
with $k$ space directions and $l$ time directions has at least
$k$ space-like diameters and at least $l$ time-like ones, see Section \ref{diamsect}.

In Section \ref{flow} we introduce pseudo-Euclidean billiards. They can be regarded 
as a particular case of projective billiards introduced in \cite{Ta1}.
The corresponding pseudo-Euclidean billiard 
map has a variational origin and exhibits  peculiar properties.
For instance, there are special (``singular'') points, where the normal to
the reflecting surface is tangent to the surface itself 
(the phenomenon impossible for Euclidean reflectors), at which the billiard
map is not defined. These points can be of two different types, and the 
reflection near them is somewhat similar to the reflection in the two different
wedges, with angles $\pi/2$ and $3\pi/2$,  in a Euclidean space.
As an illustration, we study in detail the case of a circle on a Lorentz plane, see Section \ref{case}.

We prove a Lorentz  version of the  Clairaut theorem on the complete integrability of the geodesic flow on a surface of revolution.
We also prove  pseudo-Euclidean analogs of the 
Jacobi-Chasles theorems and show the integrability of 
the billiard in the ellipsoid and the geodesic flow on the ellipsoid
in a pseudo-Euclidean space. Unlike the Euclidean situation, 
the number of ``pseudo-confocal" conics passing through
a point in pseudo-Euclidean space may differ for different points 
of space, see Section \ref{psellipt}.

Throughout the paper, we mostly refer
to  ``pseudo-Euclidean spaces'' or ``pseudo-Riemannian manifolds'' to emphasize
arbitrariness of the number of space- or time-like directions. ``Lorentz" means 
that the signature is of the form $(k,1)$ or $(1,l)$.
Note that the contact structure on null geodesics 
was previously known, at least, for the Lorentz case -- see \cite{Lo, N-T}, and it had been important in various  causality questions in the physics literature.
Apparently, pseudo-Euclidean billiards have not been  considered before, 
nor was the integrability of pseudo-Riemannian geodesic flows on 
quadratic surfaces different from pseudospheres.

\bigskip

{\bf Acknowledgments}. It is our great pleasure to thank Max-Planck-Institut in Bonn for its support and hospitality. We are grateful to  J. C. Alvarez, V. Arnold, D. Genin, C. Duval , P. Lee and especially  P. Iglesias-Zemmour for stimulating discussions and to D. Genin for numerical study of Lorentz billiards. The first author was partially supported by an NSERC grant and the second by an NSF grant. This research was partially conducted
during the period the first author was employed by the Clay Mathematics Institute as a 
Clay Book Fellow.


\section{Symplectic and contact structures on the spaces of oriented geodesics} \label{str}

\subsection{General construction} \label{gen}
Let $M^n$ be a smooth manifold with a pseudo-Riemannian metric $\langle\, ,\rangle$ of signature $(k,l),\ k+l=n$. Identify the tangent and cotangent spaces via the metric. Let $H: T^*M\to \R$ be the Hamiltonian of the metric: $H(q,p)=\langle p,p\rangle/2$. The geodesic flow in $T^*M$ is the Hamiltonian vector field $X_H$ of $H$.  

A geodesic curve in $M$ is a projection of a trajectory of $X_H$ to $M$. Let ${\cal L}_+, {\cal L}_-, {\cal L}_0$ be the spaces of oriented non-parameterized space-, time- and light-like geodesics (that is, $H=const>0, <0$ or $=0$, respectively). Let ${\cal L}={\cal L}_+\cup   {\cal L}_-\cup  {\cal L}_0$ be the space of all oriented geodesic lines.
We assume that these spaces are smooth manifolds (locally, this is always the case); then ${\cal L}_0$ is the common boundary of ${\cal L}_\pm$.

Consider the actions of $\R^*$ on the tangent and cotangent bundles by rescaling (co)vectors. The Hamiltonian $H$ is homogeneous of degree 2 in the variable $p$.  Refer to this action as the dilations. Let $E$ be the Euler field in $T^*M$ that generates the dilations.

\begin{theorem} \label{sympcont}
The manifolds ${\cal L}_\pm$ carry symplectic structures obtained from $T^*M$ by Hamiltonian reduction on the level hypersurfaces $H=\pm 1$. The manifold ${\cal L}_0$ carries a contact structure whose symplectization is the Hamiltonian reduction of the symplectic structure in $T^*M$ (without the zero section) 
on the level hypersurface $H=0$. 
\end{theorem}

\proof Consider three level hypersurfaces: $N_{-1}=\{H=-1\}, N_0=\{H=0\}$ and $N_{1}=
\{H=1\}$. The Hamiltonian reduction on the first and the third yields the symplectic structures in ${\cal L}_\pm$. This is the same as in the Riemannian case, see, e.g., \cite{A-G}.

Consider the hypersurface $N_0$ in $T^*M$ without the zero section. We have two vector fields on it, $X_H$ and $E$, satisfying $[E,X_H]=X_H$. Denote the Hamiltonian reduction of $N_0$ by $P$, it is the quotient of $N_0$ by the $\R$-action with the generator $X_H$ (sometimes, $P$ is called the space of {\it scaled} light-like geodesics). Then  ${\cal L}_0$ is the  quotient of $P$ by the dilations; denote the projection $P \to {\cal L}_0$ by $\pi$. Note that $E$ descends on $P$ as a vector field $\bar E$. Denote by $\bar \omega$ the symplectic form on $P$. Let $\bar \lambda=i_{\bar E} \bar \omega$. We have:
$$
d \bar \lambda=\bar \omega,\  L_{\bar E} (\bar \omega)=\bar \omega, \ 
L_{\bar E} (\bar \lambda)=\bar \lambda.
$$
Thus $(P,\bar \omega)$ is a homogeneous symplectic manifold with respect to the Euler field $\bar E$. Consider the distribution Ker $\bar \lambda$ on $P$. Since $\bar E$ is tangent to this distribution, Ker $\bar \lambda$ descends to a distribution on ${\cal L}_0$. This is
a contact structure whose symplectization is $(P,\bar \omega)$.

To prove that the distribution on ${\cal L}_0$ is indeed contact, let $\eta$ be a local 1-form defining the distribution. Then $\pi^*(\eta)=\bar \lambda$. Hence  
$$
\pi^*(\eta\wedge d\eta^{n-2})=\bar \lambda\wedge \bar \omega^{n-2}=\frac{1}{n-1}\ 
i_{\bar E} \bar \omega^{n-1}.
$$
Since $\bar \omega^{n-1}$ is a volume form, the last form does not vanish.
\proofend


\subsection{Examples} \label{subex}
\begin{example} \label{plane}
{\rm Let us compute the area form on the space of lines in the Lorentz plane
with the metric $ds^2=dx dy$. A vector $(a,b)$ is orthogonal to $(a,-b)$. Let $D(a,b)=(b,a)$ be the linear operator identifying vectors and covectors via the metric.
  
The light-like lines are horizontal or vertical, the space-like have positive and the time-like negative slopes. Each space ${\cal L}_+$ and ${\cal L}_-$ has two components. To fix ideas, consider space-like lines having the direction in the first coordinate quadrant.
 Write the unit directing vector of a line  as $(e^{-u},e^u),\ u\in \R$.  Drop the perpendicular $r(e^{-u},-e^u),\ r\in \R$, to the line from the origin. Then $(u,r)$ are coordinates in ${\cal L}_+$. Similarly one introduces coordinates in ${\cal L}_-$.
 
 \begin{lemma} (cf.\cite{Bir,dO-S}) \label{symstr2}
 The area form $\omega$ on  ${\cal L}_+$ is equal to $2 du\wedge dr$, and to $-2 du\wedge dr$ on ${\cal L}_-$.
 \end{lemma}

\proof
Assign to a line with coordinates $(u,r)$ the covector $p=D(e^{-u},e^u)=(e^u,e^{-u})$ and the point $q=r(e^{-u},-e^u)$. This gives a section of the bundle $N_1\to {\cal L}_+$, and the symplectic form $\omega=dp\wedge dq$ equals, in the $(u,r)$-coordinates, $2 du\wedge dr$. The computation for ${\cal L}_-$ is similar.
\proofend
}
\end{example}

\begin{example} \label{desit1}
{\rm Consider the  Lorentz space with the metric $d x^2 + d y^2 - d z^2$; let $H^2$ be the upper sheet of the hyperboloid $ x^2 +y^2 -z^2=-1$ and $H^{1,1}$ the hyperboloid of one sheet $x^2+y^2-z^2 =1$. The restriction of the ambient metric to $H^2$ gives it a Riemannian metric of constant negative curvature and the restriction to $H^{1,1}$ a Lorentz metric of constant curvature. The geodesics of these metrics are the intersections of the surfaces with the planes through the origin; the light-like lines of $H^{1,1}$ are the straight rulings of the hyperboloid. The central projection on a plane induces a (pseudo)-Riemannian metric therein whose geodesics are straight lines (for $H^2$, this is the Beltrami-Klein model of the hyperbolic plane).

The scalar product in the ambient space determines duality between lines and points by assigning to a vector the orthogonal plane. In particular, to a point of $H^2$ there corresponds a space-like line in $H^{1,1}$ (which is a closed curve).
More precisely, $H^2$ (which is the upper sheet of the hyperboloid)
is identified with the space of positively (or ``counterclockwise'') oriented 
space-like lines of $H^{1,1}$, while the lower sheet of the same hyperboloid (that is, $H^2$ with the opposite orientation) is identified with the space of 
negatively  oriented  lines. On the other hand, $H^{1,1}$ is identified with the space of oriented lines in $H^2$. The space of oriented time-like lines of $H^{1,1}$ (which are not closed) is also identified with $H^{1,1}$ itself. The area forms on the spaces of oriented lines coincide with the area forms on the respective surfaces, induced by the ambient metric. 

This construction is   analogous to the projective duality between points and oriented great circles of the unit sphere in 3-space.
 }
\end{example}


\subsection{Pseudo-Euclidean space} \label{psE} 
Let $V^{n+1}$ be a vector space with an indefinite non-degenerate quadratic form. Decompose $V$ into the orthogonal  sum of the positive and negative subspaces; denote by $v_1,v_2$ the positive and negative components of a vector $v$, and likewise, for covectors. The scalar product in $V$ is given by the formula 
$\langle u,v\rangle=u_1\cdot v_1 - u_2\cdot v_2$ where $\cdot$ is the Euclidean dot product. Let $S_{\pm}$ be the unit pseudospheres in $V$ given by the equations $|q_1|^2-|q_2^2|=\pm 1$.

The next result and its proof are similar to the familiar Euclidean case.

\begin{proposition} \label{pmsymp}
 ${\cal L}_{\pm}$ is (anti)symplectomorphic to $T^* S_{\pm}$. 
 \end{proposition}

\proof Consider the case of  ${\cal L}_+$. Assign to a space-like line $\ell$ its unit vector $v$, so that  $|v_1|^2-|v_2^2|=1$, and a point $x\in \ell$ whose position vector is orthogonal to $v$, that is, $\langle x,v\rangle=0$. Then $v\in S_+$ and $x\in T_v S_+$. Let $\xi\in T^*_v S_+$ be the covector corresponding to the vector $x$ via the metric: $\xi_1=x_1,\xi_2=-x_2$. Then the canonical symplectic structure in $T^* S_+$ is $d\xi\wedge dv= dx_1\wedge dv_1-dx_2\wedge dv_2$.

The correspondence $\ell\mapsto (q,p)$, where $q=x$ and $p=(p_1,p_2)=(v_1,-v_2)$ is the covector corresponding to the vector $v$ via the metric, is a section of the bundle $N_1\to {\cal L}_+$. Thus the symplectic form $\omega$ on ${\cal L}_+$ is the pull-back of the form $dp\wedge dq$, that is,  $\omega=dv_1\wedge dx_1 - dv_2\wedge dx_2$. Up to the sign, this is the symplectic structure in $T^* S_+$.
\proofend

A light-like line is characterized by its point $x$ and a vector $v$ along the line; one has $\langle v,v\rangle=0$. The same line is determined by the pair $(x+sv,tv),\ s\in\R, t\in \R_+^*$. The respective vector fields $v\partial x$ and $v\partial v$ are the Hamiltonian and the Euler fields, in this case. 

We shall now describe the contact structure in  ${\cal L}_0$  geometrically. 

Assign to a line $\ell \in {\cal L}_0$ the set $\Delta(\ell)\subset {\cal L}$ consisting of the oriented lines in the affine hyperplane, orthogonal to $\ell$. Then $\ell \in \Delta(\ell)$ and $\Delta(\ell)$ is a smooth $(2n-2)$-dimensional manifold, the space of oriented lines in $n$-dimensional space. Denote by $\xi(\ell) \subset T_{\ell} {\cal L}$ the tangent hyperplane to $\Delta(\ell)$ at point $\ell$. 

Denote by $S_0$ the spherization of the light cone: $S_0$ consists of equivalence classes of non-zero vectors $v\in V$ with $\langle v,v\rangle=0$ and $v\sim tv, t>0$. Let $E$ be the 1-dimensional $\R_+^*$-bundle over $S_0$ whose sections are functions $f(v)$, homogeneous of degree $1$. Denote by $J^1 E$ the space of 1-jets of sections of $E$; this is a contact manifold.

\begin{proposition} \label{zerocont}
1. $\xi(\ell)$ is the contact hyperplane of the contact structure in ${\cal L}_0$.\\
2. ${\cal L}_0$ is contactomorphic to $J^1 E$.\footnote{In the case of a Lorentz space, ${\cal L}_0$ is also contactomorphic to the space of cooriented contact elements of a Cauchy surface, see \cite{Lo}.}
\end{proposition}

\proof By construction of Theorem  \ref{sympcont}, the contact hyperplane at $\ell$ is the projection to $T_{\ell} {\cal L}_0$ of the kernel of the Liouville form $vdx$ (identifying vectors and covectors via the metric). Write an infinitesimal deformation of $\ell=(x,v)$ as
$(x+\eps y,v+\eps u)$. This is in Ker $vdx$ if and only if $\langle y,v\rangle=0$. The deformed line is light-like, hence $\langle v+\eps u,v+\eps u\rangle=0$ mod $\eps^2$, that is, $\langle u,v\rangle=0$. Thus both the foot point and the directional vector of the line $\ell$ move in the hyperplane, orthogonal to $\ell$, and therefore the contact hyperplane at $\ell$ is contained in $\xi(\ell)$. Since the dimensions coincide, $\xi(\ell)$ is this contact hyperplane. In particular, we see that $\Delta(\ell)$ is tangent to ${\cal L}_0$ at point $\ell$. This proves the first statement.

Assign to $\ell=(x,v)$ the 1-jet of the function $\phi(\ell)=\langle x,\cdot\rangle$ on $S_0$. This function is homogeneous of degree 1. The function $\phi(\ell)$ is well defined: since $v$ is orthogonal to $v$ and to $T_v S_0$, the function $\phi(\ell)$ does not change if $x$ is replaced by $x+sv$. Thus we obtain a diffeomorphism $\phi: {\cal L}_0 \to J^1 E$.

To prove that $\phi$ preserves the contact structures, let $f$ be a test section of $E$. By definition of the contact structure in $J^1 E$, the 1-jet extension of $f$ is a Legendrian  manifold. Set: $x(v)=\nabla f(v)$ (gradient taken with respect to the pseudo-Euclidean structure). We claim that $\phi(x(v),v)=j^1 f (v)$. Indeed, by the Euler formula, 
\begin{equation} \label{eu}
\langle x(v),v\rangle=\langle \nabla f(v),v\rangle=f(v),
\end{equation}
that is,  the value of the function $\langle x(v),\cdot\rangle$ at point $v$ is $f(v)$. Likewise, let $u\in T_v S_0$ be a test vector. Then the value of the differential $d\langle x(v),\cdot\rangle$ on $u$ is $\langle \nabla f(v),u\rangle=df_v(u)$.

It remains to show that the manifold $\phi^{-1}(j^1f)=\{(x(v),v)\}$ is Legendrian in ${\cal L}_0$. Indeed, the contact form is $v dx$. One has:
$$
v d(x(v))=d\langle x(v),v\rangle -x(v)dv=df-\nabla f dv =0;
$$
the second equality is due to (\ref{eu}).  Therefore $\phi^{-1}(j^1f)$ is a Legendrian submanifold, and the second claim follows.
\proofend


\subsection{Symplectic, Poisson and contact structures} \label{Pois}

The contact manifold ${\cal L}_0$
is the common boundary of the two open symplectic manifolds ${\cal L}_\pm$.
Suppose that $n\ge 2$, that is,  we consider lines
in at least three-dimensional space $V^{n+1}$.

\begin{theorem}\label{poisson}
Neither the symplectic structures of  ${\cal L}_\pm$, nor their inverse
Poisson structures, extend smoothly across  the boundary ${\cal L}_0$
to the corresponding structure on the total  space 
${\cal L}={\cal L}_+\cup  {\cal L}_0\cup   {\cal L}_-$.
\end{theorem}

\begin{remark}
{\rm When $n=1$ the symplectic structures go to infinity
as we approach the one-dimensional manifold ${\cal L}_0$. The corresponding
Poisson structures, which are inverses of the symplectic ones, extend
smoothly across  ${\cal L}_0$.

This can be observed already in 
the explicit computations of Example  \ref{subex}. Recall that for the metric
$ds^2=dxdy$ in $V^{2}$ and the lines directed by vectors $(e^{-u},e^u), \,\,u\in\R$
the symplectic structure in the corresponding coordinates $(u,r)$ at 
${\cal L}_+$ has the form $2du\wedge dr$, see Lemma \ref{symstr2}.

Now consider a neighborhood of a light-like line among all lines, that is,
a neighborhood of a point in ${\cal L}_0$ regarded as a boundary submanifold 
between ${\cal L}_+$ and ${\cal L}_-$. Look at the variation
$\xi_\eps=(1,\eps)$ of the horizontal (light-like) direction
$\xi_0=(1,0)$, and regard $(\eps, r)$ as the coordinates in this neighborhood.
For $\eps>0$ the corresponding half-neighborhood lies in ${\cal L}_+$,
while the coordinates $u$ and $\eps$ in this half-neighborhood are related
as follows. Equating the slope of $(1,\eps)$ to the slope of
$(e^{-u},e^u)$ we obtain the relation $\eps=e^{2u}$ 
or $u=\frac 12\ln\eps$. Then the symplectic structure
$\omega =2du\wedge dr=d\ln\eps\wedge dr=\frac 1\eps\, d\eps\wedge dr$.
One sees that $\omega\to\infty$ as $\eps\to 0 $.
The Poisson structure, inverse to $\omega$, is given by the bivector
field $\eps\frac{\partial}{\partial\eps}\wedge
\frac{\partial}{\partial r}$ and it extends smoothly across the border 
$\eps=0$.}
\end{remark}

\begin{example}
{\rm 
Let us compute the symplectic strictures on lines in the 3-dimensional space $V^{3}$
with the metric $dxdy-dz^2$. We parametrize the space-like directions
by $\xi=(e^{-u}\cosh\phi,e^{u}\cosh\phi, \sinh \phi)$, where
$u\in\R,\,\phi\in\R$.
The operator $D$ identifying vectors and covectors has the form $D(a,b,c)= 
(b/2,a/2,-c)$. Choose the basis of vectors orthogonal to $\xi$ as
$$e_1=(e^{-u}\sinh\phi,e^{u}\sinh\phi, \cosh \phi)\,\,\text{ and }\,\,
e_2=(e^{-u},-e^{u}, 0)\,.
$$
The symplectic structure $\omega=dp\wedge dq$ for $q=r_1e_1+r_2e_2$ and 
$p=D\xi=(e^{u}\cosh\phi/2, e^{-u}\cosh\phi/2,-\sinh \phi)$  has the following explicit expression in coordinates
$(u,\phi,r_1,r_2)$:
$$
\omega=-\,d\phi\wedge dr_1
+\cosh\phi \,du\wedge dr_2
-r_2\sinh\phi \,d\phi\wedge du\,.
$$
}
\end{example}

Now we are in a position to prove  Theorem \ref{poisson} on non-extendability.

\proof
The impossibility of extensions follows from
the fact that the ``eigenvalues'' of the symplectic structures $\omega$
of  ${\cal L}_\pm$ go to both 0 and $\infty$, as we approach ${\cal L}_0$ 
from either side. (Of course, according to the Darboux
theorem, the eigenvalues of the symplectic structures are not well defined,
but their zero or infinite limits are.) More precisely, let
$\alpha=\sum a_{ij}dx_i\wedge dx_j$ be a meromorphic
2-form written in local coordinates
$\{x_i\}$ in a neighborhood of a point $P$.

\begin{lemma} \label{eigen1}
The number of eigenvalues of the matrix $A=(a_{ij})$
which go to 0 or $\infty$ as $x\to P$ does not depend on the choice
of coordinates $\{x_i\}$.
\end{lemma}

\proof
Indeed, under a coordinate change $x=\eta(y)$, the matrix $A$
changes to $(J\eta)^*A(J\eta)$ in coordinates $\{y_j\}$,
where $J\eta$ is the Jacobi matrix of the diffeomorphism $\eta$.
Since $J\eta$ is bounded and non-degenerate, 
this change preserves (in)finiteness or vanishing
of the limits of the eigenvalues of $A$.
\proofend

Now the theorem follows from 

\begin{lemma} \label{eigen2}
The eigenvalues of the 2-form $\omega$ in coordinates $(u,\phi,r_1,r_2)$
go to both $0$ and $\infty$ as $r_2\to\infty$ (while keeping other 
coordinates fixed).
\end{lemma}

\proof
Indeed, the matrix of $\omega$ has the following (biquadratic)
characteristic equation:
$\lambda^4 +a\lambda^2+b=0$, where $a=1+ r_2^2 \sinh^2\phi +\cosh^2\phi,\,
b=\cosh^2\phi$. As $r_2\to\infty$, so does $a$, whereas $b$ does not change. Thus the sum of the squares of the roots goes to infinity, whereas their product is constant. Hence 
the equation has one pair of roots going to 0, while the other goes to infinity.
\proofend

The limit $r_2\to\infty$ means that one is approaching the boundary
of the space ${\cal L}_+$. The infinite limit of the eigenvalues means that the
symplectic structure $\omega$ does not extend smoothly  across ${\cal L}_0$,
while the zero limit of them means that the
 Poisson structure inverse to $\omega$ is non-extendable as well. 
The case of higher dimensions $n$ can be treated similarly.
\proofend

\begin{remark}
{\rm
The contact planes in ${\cal L}_0$ can be viewed as the subspaces of directions
in the tangent spaces $T_*{\cal L}_0$, on which the limits of 
the ${\cal L}_\pm$-symplectic structures are finite.
One can also see that the existence of  extensions of
the symplectic or Poisson structures would mean the presence
of other intrinsic structures, different from the contact one, on  the boundary 
${\cal L}_0$.
Indeed,  the existence of a symplectic structure extension would imply 
the existence of a presymplectic structure (and hence, generically, 
a characteristic direction field), rather than of a contact distribution,  
on  ${\cal L}_0$. 

On the other hand, consider the Poisson structures
on ${\cal L}_\pm$ which are inverses of the corresponding symplectic structures.
The assumption of a smooth extension of such Poisson structures would 
mean the existence of a Poisson structure on ${\cal L}_0$ as well. 
The corresponding foliation of ${\cal L}_0$ by symplectic leaves would 
be integrable, while the contact distribution  is not.
}
\end{remark}


\subsection{Local Lie algebra of geodesics and its Poissonization}

It turns out that the space of all pseudo-Riemannian geodesics 
(i.e., including all three types: space-, time-, and light-like 
ones) has a  structure of a Jacobi manifold, or a local Lie algebra.
This structure is not canonical and it depends on the choices described below, but it shows how  symplectic (${\cal L}_\pm$) and contact  ($ {\cal L}_0$) manifolds
constituting $\cal L$ smoothly fit together in the framework of a Jacobi manifold. 

Recall that a manifold is said to have a Jacobi structure   if the space of functions on it 
(or, more generally, the space of sections of a line bundle over it) is
equpped with a Lie bracket (a bilinear skew-symmetric operation
satisfying the Jacobi identity) which is local over the manifold.
Locality of the bracket means that it is defined by differential operators
on functions or sections \cite{Kir}.
A.~Kirillov proved that such a manifold
naturally decomposes into a union of presymplectic and contact manifolds
with natural ``pre-Poisson" or Lagrange brackets on them, respectively, 
\cite{Kir}.
A presymplectic manifold is a manifold equipped with a 2-form which is 
the product of a symplectic 2-form and a nonvanishing function.

In \cite{DLM} it was shown that any Jacobi manifold can be obtained from a homogeneous 
Poisson manifold of dimension one  bigger, called its Poissonization,
by choosing a hypersurface in the latter. It turns out that the space 
${\cal L}$ of all geodesics on a pseudo-Riemannian manifold $M$ has a natural Poissonization with a simple canonical structure. As in Section \ref{str}, we assume that the spaces of geodesics 
are smooth manifolds (or consider the local situation).

\begin{theorem}\label{locallie}
The space ${\cal L}={\cal L}_+\cup  {\cal L}_0\cup   {\cal L}_-$
of all geodesics on a pseudo-Riemannian manifold $M$ can be identified with the
quotient of a homogeneous Poisson manifold with respect to dilations.
The images of the symplectic leaves of this Poisson manifold in the quotient 
correspond to the spaces ${\cal L}_+,  {\cal L}_0,$ and ${\cal L}_-$.
\end{theorem}

\proof 
Consider $T^*M$ (without the zero section) with the standard 
symplectic structure. Abusing the notation, in the proof below we denote it
by the same symbol $T^*M$, and use other 
notations of Theorem \ref{sympcont}.
Let $H=\langle p,p\rangle/2$ be the Hamiltonian, and $X_H$ 
the corresponding Hamiltonian field.

Now we consider the manifold $PM:= T^*M/X_H$, 
that is, instead of first 
confining ourselves to the levels of $H$ as in  Theorem \ref{sympcont}, 
we  take the quotient with respect to the $\R$-action of 
the Hamiltonian field $X_H$ right away. 
Then  $PM$ is a Poisson manifold
as a quotient of the Poisson manifold  $T^*M$ along the action of the field 
$X_H$, which respects the Poisson structure. 
Furthermore, the symplectic leaves of $PM$ 
have codimension 1, since  $PM$ was obtained as a quotient of a symplectic
manifold (i.e., nondegenerate a Poisson  structure) 
by a one-dimensional group. These leaves
are exactly the levels of $H$ in $PM$.

We claim that the manifold $PM$ can be regarded as the Poissonization 
of the space ${\cal L}$ of all geodesics: it has a natural  
Poisson structure, homogeneous with respect to the action of dilations, 
and such that the quotient space coincides with ${\cal L}$: ${\cal L}=PM/\R^*$.
Indeed, consider the action of $\R^*$-dilations $E$ on $PM$. 
It is well defined on $PM$ due to the relation $[E,X_H]=X_H$. 
Note that the symplectic leaves, i.e. $H$-levels, 
are transversal to $E$ wherever $H\not=0$, and are tangent to $E$ when $H=0$.

For  $H\not=0$, the quotient space with respect to the  $\R^*$-action $E$
can be described by the levels $H=\pm1$, and the latter
correspond to the spaces of space-like or time-like geodesics.
(Note that here we have made the same Hamiltonian reduction
as in Theorem \ref{sympcont}, but in the oposite 
order: first taking the quotient, and then passing to the restriction.)

For $H=0$, we have one leaf with the field $E$ in it, which exactly constitutes
the setting for defining the space ${\cal L}_0$ of light-like
geodesics in the proof of Theorem \ref{sympcont}, 
where this leaf is called $P$, and the field is $\bar E$.
This leads to the contact structure on ${\cal L}_0$ after taking the quotient.
\proofend

\begin{corollary}
The  quotient space ${\cal L}=PM/\R^*$ can be (locally) 
endowed with a smooth Jacobi structure upon choosing any section
of this $\R^*$-bundle over ${\cal L}$.
\end{corollary}

\begin{remark}
{\rm The formulas relating the homogeneous Poisson structures and the Jacobi structures 
on the sections in the general setting can be found in \cite{DLM}.
The smooth Jacobi structure on the manifold ${\cal L}={\cal L}_+\cup  {\cal L}_0\cup   {\cal L}_-$
depends on the choice of the hypersurface in $PM$ realizing the section.
To describe the structure ``independent of this choice,'' one
can consider  the Poisson structure on ${\cal L}_\pm$, 
which, up to conformal changes, gives rise to a 
``conformal cosymplectic structure" and
captures many features of the neighboring symplectic structures.
This approach is developed in \cite{Patrick}.

On the other hand,  as we have seen in the preceding section, 
the symplectic structures of the two open submanifolds
${\cal L}_\pm$ blow up as one approaches
their common contact boundary ${\cal L}_0$. 
To see how the coexistence of a smooth Jacobi structure
and the blowing up symplectic structures fit together, we consider 
the corresponding homogeneous cone $P=\{H=0\}$ 
(consisting of scaled null geodesics) 
over the space of null geodesics.  Then the space 
${\cal L}_+$ of space--like geodesics ``approaches the cone $P$ 
at infinity'' (same for ${\cal L}_-$). 
Now consider the family of spaces $\{H=\lambda\}$, 
all isomorphic to the set of space-like geodesics. 
The picture is similar to a family of 
hyperboloids approaching the quadratic cone.  At any given point of the
cone the convergence will be smooth, since we are taking different (closer
and closer) hyperboloids, while the corresponding structures on the
quotients, when we fix one hyperboloid, and where the structure ``comes
from infinitely remote points" will not necessarily have a nice convergence.

One can show that this type of blow-up of the symplectic structures
in Jacobi manifolds is typical: there is a version of the Darboux theorem
showing that locally all such degenerations look alike, cf. \cite{MZ}.
}
\end{remark}

\subsection{Hypersurfaces and submanifolds} \label{hyps}

Let $M\subset V$ be an oriented smooth hypersurface. Assign to a point $x\in M$ its oriented normal line $\ell(x)$. We obtain a Gauss map $\psi: M\to {\cal L}={\cal L}_+\cup   {\cal L}_-\cup  {\cal L}_0$. Denote by $\psi_+,\psi_-$ and $\psi_0$ its space-, time- and light-like components, i.e. the restriction of the 
map $\psi$ to those parts of $M$, where the normal  $\ell(x)$ is respectively space-, time- or light-like. Note that if the normal line is light-like then it is tangent to the hypersurface.

\begin{proposition} \label{gauss}
The images of $\psi_{\pm}$ are Lagrangian and the image of $\psi_0$ is Legendrian.
\end{proposition}

\proof Consider  the case of $\psi_+$ (the time-like case being similar). Denote by $\nu(x)$ the ``unit" normal vector to $M$ at point $x$ satisfying $\langle \nu(x),\nu(x)\rangle=1$. The line $\ell(x)$ is characterized by its vector $\nu(x)$ and its point $x$; the correspondence $\ell\mapsto (\nu(x),x)$ is a section of the bundle $N_1\to {\cal L}_+$ over the image of $\psi_+$. We need to prove that the form $dx\wedge d\nu(x)$ vanishes on this image. Indeed, the 1-form $\nu(x) dx$ vanishes on $M$ since $\nu(x)$ is a normal vector, hence its differential is zero as well.

In the case of $\psi_0$, we do not normalize $\nu(x)$. The 1-form $\nu(x) dx$ is still zero on $M$, and this implies that the image of the Gauss map in ${\cal L}_0$ is Legendrian, as in the proof of Proposition \ref{zerocont}.
\proofend

\begin{remark} \label{imm}
{\rm The maps $\psi_{\pm}$ are immersions but $\psi_0$ does not have to be one. For example, let $M$ be a  hyperplane such that the restriction of the metric to $M$ has a 1-dimensional kernel. This kernel is the normal direction to $M$ at each point.  These normal lines foliate $M$, the leaves of this foliation are the fibers of the Gauss map $\psi_0$, and its image is an $n-1$-dimensional space.
}
\end{remark}

\begin{remark} 
{\rm More generally, let $M\subset V$ be a smooth submanifold of any codimension.
Assign to a point $x\in M$ the set of all oriented normal lines to $M$ at $x$.
This also gives us  a Gauss map $\psi: {\cal N}_M\to {\cal L}$ of the 
normal bundle ${\cal N}_M$ of the submanifold $M\subset V$ into  ${\cal L}$
with space-, time- and light-like components $\psi_+,\psi_-$ and $\psi_0$ 
respectively. In this setting, Proposition \ref{gauss} still holds,  while
the proof requires only cosmetic changes.

Note that the Jacobi structure approach, discussed in the last section,
explains why one obtains Lagrangian~/~Legendrian submanifolds by considering 
spaces of normals to various
varieties in a pseudo-Riemannian space: they are always Lagrangian submanifolds
in the Poissonization, before taking the quotient.
}
\end{remark}

\begin{example} \label{fiber}
{\rm The set of oriented lines through a point provides an example of a submanifold in ${\cal L}$ whose intersection with ${\cal L}_+ \cup {\cal L}_-$ is Lagrangian and with ${\cal L}_0$ Legendrian.
}
\end{example}

\begin{figure}[hbtp]
\centering
\includegraphics[width=2in]{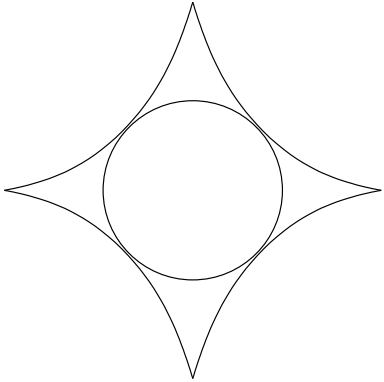}
\caption{The caustic of a circle in the Lorentz plane}
\label{astr}
\end{figure}

\begin{example} \label{cau}
{\rm Consider the circle $x^2+y^2=1$ on the Lorentz plane with the metric $dx^2-dy^2$. Then the caustic, that is, the envelope of the normal lines to the circle, is the astroid $x^{2/3}+y^{2/3}=2^{2/3}$, see figure \ref{astr} (note that the caustic of an ellipse in the Euclidean plane is an astroid too).   The role of Euclidean circles is played by the pseudocircles, the hyperbolas $(x-a)^2-(y-b)^2=c$: their caustics degenerate to points. 
}
\end{example}

\section{Billiard flow  and billiard transformation} \label{flow}

\subsection{Definition of the billiard map} \label{defmp}
Let $M$ be a pseudo-Riemannian manifold  with a smooth boundary $S=\partial M$. The billiard flow in $M$ is a continuous time dynamical system  in $TM$. The motion of  tangent vectors in the interior of $M$ is free, that is, coincides with the geodesic flow. Suppose that a vector hits the boundary at point $x$. Let $\nu(x)$ be the normal to $T_xS$. If $x$ is a singular point, that is, the restriction of the metric on $S$ is singular or, equivalently, $\langle \nu(x),\nu(x)\rangle=0$, then the billiard trajectory stops there. Otherwise the billiard reflection occurs.

Since $x$ is a non-singular point, $\nu(x)$ is transverse to $T_xS$. Let $w$ be  the velocity of the incoming point. Decompose it  into the the tangential and normal components, $w=t+n$.
Define the billiard reflection by setting $w_1=t-n$ to be the velocity of the outgoing point. Clearly $|w|^2=|w_1|^2$. In particular, the billiard reflection does not change the type of a geodesic: time-, space- or light-like.
 
We view the billiard map $T$ as acting on oriented geodesics and sending an incoming ray to the outgoing one. 

\begin{example} \label{hyp}
{\rm Let the pseudocircle $x^2-y^2=c$ be a billiard curve (or an ideally reflecting mirror) in the Lorentz plane with the metric $dx^2-dy^2$. Then all normals to this curve pass through the origin, and so every billiard trajectory from the origin reflects back to the origin. The same holds in multi-dimensional pseudo-Euclidean spaces.
}
\end{example}

\begin{example} \label{desit2}
{\rm In the framework of Example \ref{desit1}, consider two billiards, inner and outer, in the hyperbolic plane $H^2$. (The latter is an area preserving mapping of the exterior of a strictly convex curve $\gamma$ defined as follows: given a point $x$ outside of $\gamma$, draw a support line to $\gamma$ and reflect $x$ in the support point; see \cite{Ta5,Ta2}.) The duality between $H^2$ and $H^{1,1}$ transforms the inner and outer billiard systems in $H^2$ to the outer and inner billiard systems in  $H^{1,1}$. 
Given a convex closed curve in $H^2$, the dual curve in $H^{1,1}$ (consisting of the points, dual to the tangent lines of the original curve) is space-like. Thus  any outer billiard in $H^2$ provides an example of a billiard in $H^{1,1}$ whose boundary is a space-like curve.
 }
\end{example}

\begin{remark} \label{varia}
{\rm Similarly to the Riemannian case, the origin of the billiard 
reflection law is variational. One can show that 
a billiard trajectory from a fixed point $A$ to a fixed point $B$ in $M$ with reflection at point $X\in S$ is an extremal of the following variational problem:
$$
I_{\tau}(\gamma_1,\gamma_2)=\int_0^\tau \langle \gamma_1'(t),\gamma_1'(t)\rangle\ dt + \int_\tau^1 \langle \gamma_2'(t),\gamma_2'(t)\rangle\ dt
$$
where $\gamma_1(t),\ 0\leq t\leq \tau$ is a path from $A$ to a point $X$ of $S$ and $\gamma_2(t),\ \tau\leq t\leq 1$ is a path from the point $\gamma_1(\tau)$ to $B$, and where $\tau\in[0,1]$ is also a variable. 
}
\end{remark}


\subsection{Reflection near a singular point} \label{singr}
\begin{figure}[hbtp]
\centering
\includegraphics[width=3in]{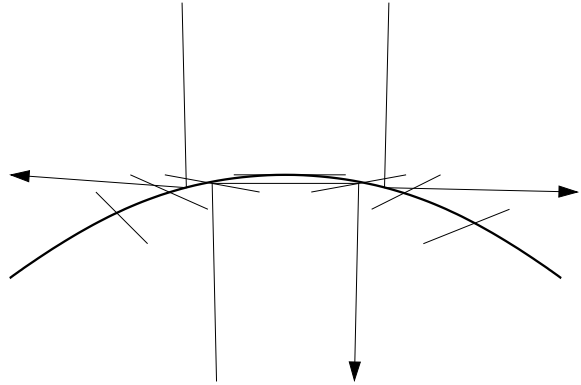}
\caption{Two types of the billiard reflection near a singular point: rays from above get scattered, while rays from the bottom have the second reflection.}
\label{norm}
\end{figure}

Let us look more carefully at the billiard reflection in a neighborhood of 
a singular point of a curve in the Lorentz plane. First of all we note that
typical singular points can be of two types, according to whether the inner
normal is oriented toward or from  the singular point as we approach it 
along the curve. These two types are shown in the same figure \ref{norm}, 
where we regard the curve as the billiard boundary either for rays coming 
from above or from below. In both cases after the reflection 
the rays are  ``squeezed'' between the tangent and the normal to the curve. 
In the first case, when  rays come from above and the normals to the curve
are oriended {\it from} the singular point, this implies that the reflected 
rays scatter away from the point. In the second case, when rays come from 
below and the normals point {\it toward} the singular point, the reflected ray
hits the boundary again.

One can see that the smooth boundary of a strictly convex domain in the Lorentz plane  has singular points of the former type  only. Indeed, up to a diffeomorphism, there exists a unique germ of normal line field at a singular point of a quadratically non-degenerate curve in the Lorentz plane -- the one shown in figure \ref{norm}. The billiard table may lie either on the convex (lower) or the concave (upper) side of the curve, whence the distinction between he two cases.

Note also that, at a singular point, the caustic of the curve always touches the curve
(cf. Example \ref{cau}). The above two cases differ by the location of the 
caustic: it can touch the curve from (a) the exterior  or (b) the interior 
of the billiard domain. The billiard inside a circle in the Lorentz plane 
has singular points of the former type only, cf. figure \ref{astr}.

The billiard reflections are drastically different in these two cases. 
In  case (b), a generic family of rays gets dispersed in opposite directions
on different sides from the singular point.
In  case (a), the situation is quite different: the scattered trajectories
are reflected toward the singular point and hit the curve one more time
in its  vicinity. Thus one considers the square of the billiard map $T^2$.

\begin{proposition} \label{singref}
Assume that a smooth billiard curve $\gamma$ in the Lorentz plane is quadratically non-degenerate at a singular point $O$. Consider a parallel beam of lines $\{\ell\}$ reflecting in an arc of $\gamma$ near point $O$,  on the convex side. Then, as the reflection points tend to $O$, the lines $T^2(\ell)$ have a limiting direction, and this direction is parallel to $\ell$.
\end{proposition}

\proof Let the metric be $dx dy$. In this metric, a vector $(a,b)$ is orthogonal to $(-a,b)$.
We may assume that the singular point is the origin, and that $\gamma$ is the graph $y=f(x)$ where $f(0)=f'(0)=0$ and $f''(0)>0$. Consider a downward incoming ray with slope $u$ reflecting in $\gamma$ at point $(s,f(s))$, then at point $(t,f(t))$, and escaping upward with slope $v$. 

\begin{figure}[hbtp]
\centering
\includegraphics[width=2in]{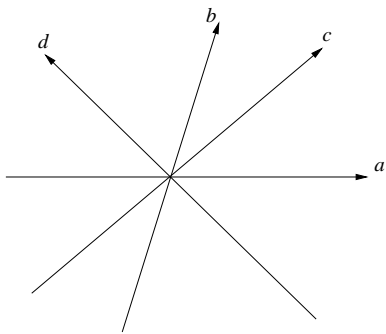}
\caption{Harmonic quadruple of lines given by four vectors}
\label{harmo1}
\end{figure}

The law of billiard reflection can be formulated as follows: the incoming billiard ray, the outgoing one, the tangent line and the normal to the boundary of the billiard table at the impact point constitute a harmonic quadruple of lines. See \cite{Ta1} for a study of projective billiards. The following criterion is convenient, see \cite{Ta1}. Let the lines be given by vectors $a,b,c,d$, see figure \ref{harmo1}. Then the lines constitute a harmonic quadruple if and only if 
\begin{equation} \label{harmcr}
[a,c][b,d]+[a,d][b,c]=0
\end{equation}
where $[,]$ is the cross product of two vectors.

For the first reflection, we have
$$
a=(1,f'(s)),\ b=(-1,f'(s)),\ c=(t-s,f(t)-f(s)),\ d=(u,1).
$$
Substitute to (\ref{harmcr}) and compute the determinants to obtain:
$$
u(f'(s))^2 (t-s)=f(t)-f(s).
$$
Similarly, for the second reflection, we have:
$$
v(f'(t))^2 (t-s)=f(t)-f(s),
$$
and hence 
\begin{equation} \label{u=v}
u(f'(s))^2=v(f'(t))^2.
\end{equation}
Write $f(x)=ax^2+O(x^3)$, then  $f'(x)=2ax+O(x^2)$,
and 
$$
\frac{f(t)-f(s)}{t-s}=a(s+t)+O(s^2,st,t^2).
$$
The above quantity equals $u(f'(s))^2$ which is $O(s^2)$, hence $t=-s+O(s^2)$. It follows that 
$(f'(t))^2=(f'(s))^2=4a^2 s^2+O(s^3)$ and, by (\ref{u=v}), that $v=u$.
\proofend

\begin{figure}[hbtp]
\centering
\includegraphics[width=2.5in]{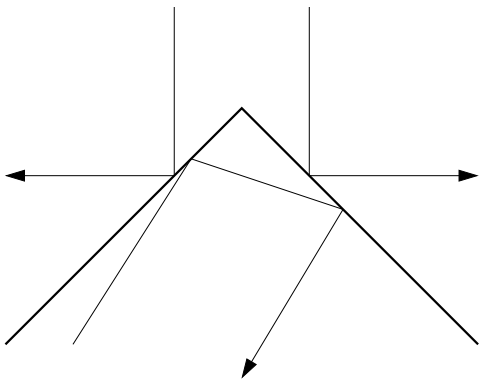}
\caption{Euclidean billiard reflection in a right angle}
\label{angle}
\end{figure}

Thus a ray meeting a curve near a singular point emerges, after two reflections, in the opposite direction. This resembles the billiard reflection in a right angle in the Euclidean plane, see figure \ref{angle}. In contrast, the reflection of a parallel beam on the concave side of a Lorentz billiard near a 
singular point resembles  the Euclidean billiard reflection from the angle $3\pi/2$ (cf. figures 2 and 4). 
Of course, this behavior excludes the existence of smooth caustics in Lorentz billiards, cf., e.g, \cite{Ta5,Ta2} for the Euclidean case.

\subsection{Symplectic and contact properties of the billiard map} \label{symprop}
Now we discuss symplectic properties of the billiard transformation.
To fix ideas, let the billiard table be geodesically 
convex.\footnote{Alternatively, one may consider the situation locally, 
in a neighborhood of an oriented line transversally intersecting $S$ 
at a non-singular point.} Denote by ${\cal L}^0$ the set of oriented 
lines that meet $S$ at non-singular points. 
The billiard map $T$ preserves the space-,   time-, and light-like parts
of ${\cal L}^0$, so we have   billiard transformations $T_+, T_-$  and $T_0$ acting on ${\cal L}_+^0, {\cal L}_-^0$ and ${\cal L}_0^0$, respectively. The (open dense) subsets  ${\cal L}_\pm^0
\subset {\cal L}_\pm$ and $ {\cal L}_0^0\subset {\cal L}_0$ carry the same 
symplectic or contact structures as the ambient spaces.

\begin{theorem} \label{invstr} 
The transformations $T_+$ and $T_-$ are symplectic and $T_0$ is a contact transformation.
\end{theorem}

\proof We adopt the approach of R. Melrose \cite{Me1,Me2}; see also \cite{A-G}. Identify tangent vectors and covectors via the metric. We denote vectors by $v$ and covectors by $p$.

Consider first the case of space-like geodesics (the case of time-like ones is similar). Let $\Sigma\subset T^*M$ be the hypersurface consisting of  vectors with foot-point on $S$. Let $Z=N_1\cap \Sigma$ and let $\Delta\subset Z$ consist of the vectors tangent to $S$.

Denote by $\nu(q)\in T^*_q M$ a conormal vector to $S$ at point  $q\in S$.
Consider the characteristics of the canonical symplectic form $\omega$ in $T^* M$ restricted to $\Sigma$.  We claim that these are the lines $(q, p+t\nu(q)), t \in \R$.

Indeed, in local Darboux coordinates, $\omega=dp \wedge dq$. The line 
$(q, p+t\nu(q))$ lies in the fiber of the cotangent bundle $T^*M$ over the point $q$ and
the vector $\xi = \nu(q) \partial /\partial p$ is tangent to this line. Then $i_{\xi}
\omega = \nu(q) dq$. This 1-form vanishes on $\Sigma$ since $\nu(q)$ is a conormal vector to $S$ at $q$.  Thus $\xi$ has the characteristic direction. Note that the quotient space by the characteristic foliation is $T^*S$.

Next we claim that the restriction of $\omega$ to $Z-\Delta$ is a symplectic form. Indeed, $Z-\Delta\subset N_1$ is transverse to the trajectories of the geodesic flow, that is, the leaves of the characteristic foliation of $N_1\subset T^* M$.

The intersections of $Z$ with the leaves of the characteristic foliation on $N_1$  determine an involution, $\tau$, which is free on $Z-\Delta\subset N_1$. If $(q,v)\in Z$ is a vector, let $q_1 \in S$ be the other intersection point of the geodesic generated by $(q,v)$ with $S$ and $v_1$ the vector  translated to point $q_1$ along the geodesic. Then $\tau(q,v)=(q_1,v_1)$.

Consider the intersections of $Z$ with the leaves of the characteristic foliation on $\Sigma$. We claim that this also determines an involution, $\sigma$, which is free on $Z-\Delta\subset \Sigma$. Indeed, let $(q,v)\in Z$, i.e., $q\in S, \langle v,v\rangle=1$. The characteristic line is $(q, v+t\nu(q))$, where $\nu(q)$ is a normal vector, and its intersection with $Z$ is given by the equation $\langle p+t\nu(q),p+t\nu(q)\rangle=1$. Since $\langle \nu(q),\nu(q)\rangle\neq 0$, this equation has two roots and we have an involution. One root is $t=0$, the other is different from 0 if $\langle v, \nu(q)\rangle\neq 0$, that is, $v$ is not tangent to $S$.

\begin{figure}[hbtp]
\centering
\includegraphics[width=3in]{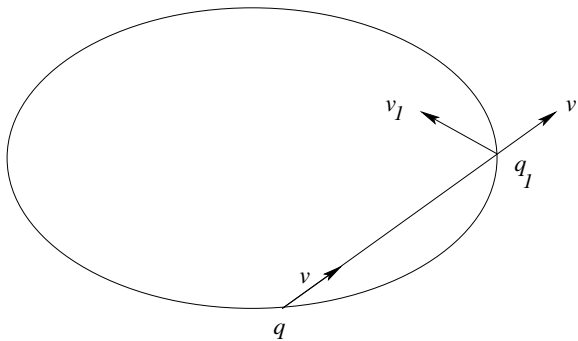}
\caption{The billiard map $F=\sigma \circ \tau$ is a composition of two involutions: $\tau(q,v)=(q_1,v), \sigma(q_1,v)=(q_1,v_1)$}
\label{2inv}
\end{figure}

Let $F=\sigma \circ \tau$; this is the billiard map on $Z$, see figure \ref{2inv}.
Since both involutions are defined by intersections with the leaves of the characteristic foliations, they preserve the symplectic structure $\omega|_Z$. Thus $F$ is a symplectic transformation of $Z-\Delta$. Let $P:Z-\Delta \to {\cal L}_+^0$ be the projection. Then $P$ is a symplectic 2-to-1 map and $P\circ F=T_+\circ P$. It follows that $T_+$ preserves the symplectic structure in ${\cal L}_+^0$.

In the case of $T_0$, we have the same picture with $N_0$ replacing $N_1$ and its symplectic reduction $P$ in place of ${\cal L}_+^0$. We obtain a symplectic  transformation of $P$ that commutes with the action of $\R_+^*$ by dilations. Therefore the map $T_0$ preserves the contact structure of ${\cal L}_0^0$.
\proofend

\begin{figure}[hbtp]
\centering
\includegraphics[width=2in]{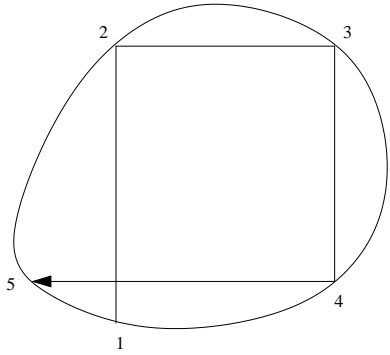}
\caption{A dynamical system on an oval}
\label{arn}
\end{figure}

\begin{remark} \label{Arn}
{\rm Consider a convex domain $D$ in the Lorentz plane with the metric $dx dy$. 
The light-like lines are either horizontal or vertical. The billiard system in $D$, restricted to light-like lines, coincides with the system described in \cite{Fest} in the context of Hilbert's 13th problem (namely, see figure \ref{arn} copied from
figure 3 on p. 8 of \cite{Fest}). The map that moves a point of the curve first  along a vertical and then along a horizontal chord is a circle map that, in case $D$ is an ellipse, is conjugated to a rotation. The same map is discussed in \cite{Sob1} in  the context of the Sobolev equation, approximately describing fluid oscillations in a fast rotating tank.
}
\end{remark}

\subsection{Diameters} \label{diamsect}
A convex hypersurface in $\R^n$ has at least $n$ diameters, which are 
2-periodic billiard trajectories in this hypersurface. 
In a pseudo-Euclidean space with signature $(k,l)$ the result is as follows.

\begin{theorem}  \label{diam}
A smooth strictly convex  closed hypersurface has at least $k$ space-like and $l$ time-like diameters.
\end{theorem}

\proof Denote the hypersurface by $Q$.
Consider the space of chords $Q\times Q$ and set $f(x,y)=\langle x-y,x-y\rangle/2$. Then $f$ is a smooth function on $Q\times Q$. The group $\Z_2$ acts on $Q\times Q$ by interchanging points, and this action is free off the diagonal $x=y$. The function $f$ is $\Z_2$-equivariant.

First we claim that a critical point of $f$ with non-zero critical value corresponds to a diameter (just as in the Euclidean case). Indeed, let $u\in T_xQ, v\in T_yQ$ be test vectors. Then $d_x f(u)=\langle x-y,u\rangle$ and $d_y f(v)=\langle x-y,v\rangle$. Since these are zeros for all $u,v$, the (non-degenerate) chord $x-y$ is orthogonal to $Q$ at both end-points. Note that such a critical chord  is not light-like, due to convexity of $Q$.

Fix a sufficiently small generic $\varepsilon >0$. Let $M\subset Q\times Q$ be a submanifold with boundary given by  $f(x,y)\geq \varepsilon$. Since the boundary of $M$ is a level hypersurface of $f$, the gradient of $f$ (with respect to an auxiliary  metric) has inward direction along the boundary, and the inequalities of Morse theory apply to $M$. Since $\Z_2$ acts freely on $M$ and $f$ is $\Z_2$-equivariant, the number of critical $\Z_2$-orbits of $f$ in $M$ is not less than the sum of $\Z_2$ Betti numbers of $M/\Z_2$.

We claim that $M$ is homotopically equivalent to $S^{k-1}$ and $M/\Z_2$ to $\RP^{k-1}$. Indeed, $M$ is homotopically equivalent to the set of space-like oriented lines intersecting $Q$. Retract this set to the set of  space-like oriented lines through the origin. The latter is the spherization of the cone $|q_1|^2>|q_2|^2$, and the projection $(q_1,q_2)\mapsto q_1$  retracts it to the sphere  $S^{k-1}$.

Since the sum of $\Z_2$ Betti numbers of $\RP^{k-1}$ is $k$, we obtain at least $k$ space-like diameters. Replacing $M$ by the manifold $\{f(x,y)\leq -\varepsilon\}$ yields $l$ time-like diameters.
\proofend

\begin{problem} \label{Push}
{\rm In Euclidean geometry, the fact that a smooth closed convex hypersurface in $\R^n$ has at least $n$ diameters has a far-reaching generalization due to Pushkar' \cite{Pu}: 
a generic immersed closed manifold $M^k \to \R^n$ has at least 
$(B^2-B+kB)/2$ diameters, that is, chords that are perpendicular to $M$ at both end-points; here $B$ is the sum of the $\Z_2$-Betti numbers of $M$. It is interesting to find a pseudo-Euclidean analog of this result.}
\end{problem}

\begin{problem} \label{Birk}
{\rm Another generalization, in Euclidean geometry, concerns the least number of periodic billiard trajectories inside  a closed smooth strictly convex hypersurface. In dimension 2, the classical Birkhoff theorem asserts that, for every $n$ and every rotation number $k$, coprime with $n$, there exist at least two $n$-periodic billiard trajectories with rotation number $k$, see, e.g., \cite{Ta5,Ta2}. In higher dimensions, a similar result was obtained recently \cite{F-T1,F-T2}. It is interesting to find analogs for billiards in pseudo-Euclidean space. A possible difficulty is that variational problems in this set-up may  have no solutions: for example, not every two points on the hyperboloid of one sheet in Example \ref{desit1} are connected by a geodesic!
}
\end{problem}


\section{The geodesic flow on a quadric and the billiard inside a quadric} \label{quad}

\subsection{Geodesics and characteristics} \label{geochar}
Let us start with a general description of the geodesics on a hypersurface in a pseudo-Riemannian manifold.

Let $M^n$ be a pseudo-Riemannian manifold and $S^{n-1}\subset M$ a smooth hypersurface. The geodesic flow on $S$ is a limiting case of the billiard flow inside $S$ when the billiard trajectories become tangent to the reflecting hypersurface. Assume that $S$ is free of singular points, that is, $S$ is a pseudo-Riemannian submanifold:
the restriction of the metric in $M$ to $S$ is non-degenerate.  The infinitesimal version of the billiard reflection law gives the following characterization of geodesics:  a geodesic on $S$ is a curve $\gamma(t)$ such that $\langle \gamma'(t),\gamma'(t)\rangle$ = const and the acceleration $\gamma''(t)$ is orthogonal to $S$ at point $\gamma(t)$ (the acceleration is understood in terms of the  covariant derivative) -- see \cite{ON}. Note that the type (space-, time-, or light-like) of a geodesic curve remains the same for all $t$.

Let $Q\subset {\cal L}$ be the set of oriented geodesics tangent to $S$. Write $Q=Q_+\cup Q_- \cup Q_0$ according to the type of the geodesics. Then $Q_{\pm}$ are hypersurfaces in the symplectic manifolds ${\cal L}_{\pm}$ and $Q_0$ in the contact manifold ${\cal L}_0$.

Recall the definition of characteristics on a hypersurface $X$ in a contact manifold $Y$, see \cite{A-G}. Assume that the contact hyperplane $C$ at point $x\in X$  is not tangent to $X$; we say that $x$ is a non-singular point. Then $C\cap T_x X$ is a hyperplane in $C$. Let $\lambda$ be a contact form. Then $\omega=d\lambda$ is a symplectic form on $C$; a different choice of the contact form, $f \lambda$, gives a proportional symplectic form $f(x) \omega$ on $C$. The characteristic line at $x$ is the skew-orthogonal complement of the hyperplane $C\cap T_x X$ in $C$.

\begin{theorem} \label{geoline}
1) The characteristics of the hypersurfaces $Q_{\pm}\subset {\cal L}_{\pm}$ consist of oriented geodesics in $M$ tangent to a fixed space- or time-like geodesic on $S$. \\
2) The hypersurface $Q_0\subset {\cal L}_0$ consists of non-singular points and its characteristics  consist of oriented geodesics in $M$ tangent to a fixed light-like geodesic on $S$.  
\end{theorem}

\proof The argument is a variation on that given in the proof of Theorem \ref{invstr} (cf. \cite{A-G}), and we use the notation from this proof. In particular, we identify vectors and covectors using the metric.

Consider $Q_+$ (the case of $Q_-$ being similar). We have the submanifold $\Delta\subset N_1$ consisting of the unit space-like vectors tangent to $S$; the projection $N_1 \to {\cal L}_+$ identifies $\Delta$ with $Q_+$. 
Likewise, the projection $\Sigma \to T^*S$ makes it possible to consider
$\Delta$  as a hypersurface in $T^*S$. The characteristics of $\Delta\subset T^*S$ are the geodesics on $S$. 

We need to show that the two  characteristic directions on $\Delta$, induced by its inclusions into $ {\cal L}_+$ and into $T^*S$, coincide. We claim that the restriction of the canonical symplectic structure $\omega$ in $T^*M$ on its codimension 3 submanifold $\Delta$ has 1-dimensional kernel at every point. If this holds then both characteristic directions on $\Delta$ coincide with these kernels and therefore with each other.

The kernel of the restriction of $\omega$ on $\Delta$ is odd-dimensional. Assume its dimension is 3; then $\Delta\subset T^*M$ is a coisotropic submanifold. We will show that this is not the case. 

Let $\nu(q)$ be the normal vector to $S$ at point $q\in S$. Since $q$ is not singular,  $\nu(q)$ is transverse to $T_q S$. Thus the vector $u=\nu(q)\ \partial/\partial q$ is transverse to $\Delta\subset T^*M$. 
So is the vector $v=\nu(q)\ \partial/\partial p$.   Let $w$ be another transverse vector such that $u,v,w$ span a transverse space to $\Delta$.
Note that 
$$
\omega(v,u)=(dp\wedge dq) (\nu(q)\ \partial/\partial p,\nu(q)\ \partial/\partial q)=\langle \nu(q),\nu(q) \rangle\neq 0.
$$
Since $\omega$ is a symplectic form, $0\neq i_{u\wedge v}\ \omega^n =C\ \omega^{n-1}$ with $C\neq 0$, and   the $2n-3$ form $i_w \omega^{n-1}$ is a volume form
on $\Delta$. This contradicts to the fact that $T_{(q,p)} \Delta$ contains a 3-dimensional subspace skew-orthogonal to $T_{(q,p)} \Delta$, and the first statement of the theorem follows.

For the second statement, we replace $N_1$ by $N_0$ and ${\cal L}_+$ by the space of scaled light-like geodesics $P$. Then $P$ and $\Delta$ are acted upon by the dilations. 
Using the notation from Theorem \ref{sympcont}, $\pi(\Delta)=Q_0$. The characteristics of $\Delta\subset P$ consist of scaled oriented geodesics tangent to a fixed light-like geodesic on $S$.

To show that the points of $Q_0\subset {\cal L}_0$ are non-singular, it suffices to prove that the hypersurface $\Delta \subset P$ is not tangent to the kernel of the 1-form $\bar \lambda$. We claim that this kernel contains the vector $v$ transverse to $\Delta$. Indeed, 
$$
\bar \lambda (v)=\bar \omega(\bar E,v)=(dp\wedge dq)(p\ \partial/\partial p, \nu(q)\ \partial/\partial p)=0,
$$
and hence ker $\bar \lambda\neq T_{(q,p)} \Delta$.

Finally, the characteristics of the conical hypersurface $\Delta=\pi^{-1}(Q_0)$ in the symplectization $P$ of the contact manifold ${\cal L}_0$ project to the characteristics of $Q_0\subset {\cal L}_0$, see \cite{A-G}, and the last claim of the theorem follows.
\proofend


\subsection{Geodesics on a Lorentz surface of revolution} \label{revolution}
Geodesics on a  surface of revolution in the Euclidean space have the following 
Clairaut first integral: $r\sin\alpha=const$, where $r$ is the distance 
from a given point on the surface to the axis of revolution, and $\alpha$ 
is the angle of the geodesic at this point with the projection of the axis
to the surface. 

Here we describe an analog of the Clairaut  integral for Lorentz surfaces 
of revolution. Let $S$ be a surface in the Lorentz space $V^{3}$ with the metric 
$ds^2=dx^2+dy^2-dz^2$ obtained by a revolution of the graph of a function $f(z)$
about the $z$-axis: it is given by the equation $r=f(z)$ for $r^2=x^2+y^2$. We assume that the restriction of the ambient metric to the surface is pseudo-Riemannian.

\begin{figure}[hbtp]
\centering
\includegraphics[width=2in]{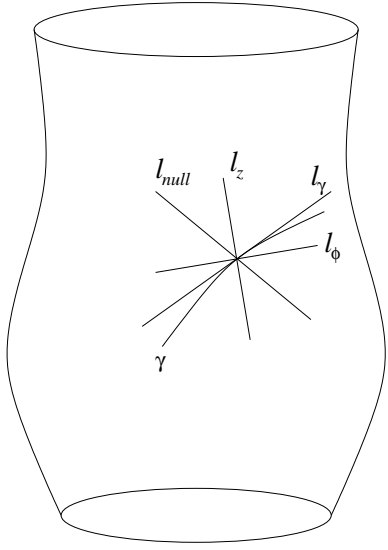}
\caption{A quadruple of tangent lines on a surface of revolution}
\label{revolv}
\end{figure}

Consider the tangent plane $T_PS$ to the surface $S$ at a point $P$ on a given 
geodesic $\gamma$. Define the following 4 lines 
in this tangent plane: the axis projection $l_{z}$ (meridian), the revolution direction
$l_\phi$ (parallel), the tangent to the geodesic $l_\gamma$, 
and one of the two null directions $l_{null}$ on the surface at the point $P$, see figure \ref{revolv}.
We denote the corresponding cross-ratio of this quadruple of lines as ${\rm cr}={\rm cr}(l_{z}, l_\phi, l_\gamma, l_{null})$.

\begin{theorem}
The function $(1-{\rm cr}^2)/r^2$ is constant along any geodesic $\gamma$
on the Lorentz surface of revolution.
\end{theorem}

\proof
The Clairaut  integral in either Euclidean or Lorentz setting
is a specification of the Noether theorem, which gives the invariance 
of the angular momenta $m=r\cdot v_\phi$ with respect to the axis of revolution.
(Here $(r,\phi)$ are the polar coordinates in the $(x,y)$-plane.)

In the Euclidean case, we have $r\phi'=|v|\sin\alpha$ and combining 
the invariance of the  magnitude of $v$ along the geodesics  
with the preservation of $m=r^2\phi'$
we immediately obtain the Clairaut  integral $r\sin\alpha=const$.

In the Lorentz setting, we first find the
cross-ratio discussed above.
Let $v$ be the velocity vector along a geodesic, and $v_r, v_\phi, v_z$ be
its radial, angle and axis projections respectively. 
Suppose the point $P$ has coordinates $(0, y_0, z_0)$;  choose 
$(x,z)$ as the coordinates in the corresponding tangent plane $T_PS$. 
The lines $l_{z}, l_\phi$, and $l_\gamma$ have 
the directions $(1,0),\,(0,1)$, and $(v_\phi, v_z)$, respectively.
The direction of null vectors in this tangent plane is the intersection
of the cone of null vectors $x^2+y^2-z^2=0$
(in the coordinates centered at $P$) with
the plane $y=f'(z_0)z$ tangent to $S$ at $P$. Thus the corresponding null 
directions are $x=\pm\sqrt{1-(f'(z_0))^2}z$.
We choose the ``plus'' direction for $l_{null}$ and find the cross-ratio
$$
{\rm cr}={\rm cr}(l_{z}, l_\phi, l_\gamma, l_{null})
=\frac{[l_{z},l_{null}][l_\gamma,l_\phi]}{[l_{z},l_\gamma][l_{null},l_\phi]}
=\frac{v_z\sqrt{1-(f'(z_0))^2}}{v_\phi}\,,
$$
that is $v_z^2(1-(f'(z_0))^2)=v_\phi^2\cdot {\rm cr}^2$.
(By choosing the other sign for the null direction $l_{null}$ we obtain 
the same relation.)

Now recall that the Lorentz length of $v$ is preserved: $v_r^2+v_\phi^2-v_z^2=1$.
Taking into account that $v_r=f'(z_0)v_z$ in the tangent plane $T_PS$, we exclude 
from this relation both $v_r$ and $v_z$ and express $v_\phi$ via the cross-ratio:
$v_\phi^2=1/(1-{\rm cr}^2)$. Thus the preservation 
of the angular momenta $m=r\cdot v_\phi$ yields the
Lorentz analog of the Clairaut theorem:
$1/m^2=(1-{\rm cr}^2)/r^2=const$ along any given geodesic on S.
\proofend

We consider the conservation of the quantity $(1-{\rm cr}^2)/r^2$, 
inverse to $m^2$, since $r\not= 0$ and this quantity makes sense even when 
${\rm cr}=\pm 1$. Note that this invariant immediately implies that a geodesic 
with a light-like initial condition stays light-like forever: along such
geodesics, ${\rm cr}^2=1$. If the geodesic is time-like ($ds^2<0$), it can be 
continued until it hits the ``tropic'' consisting of singular points.

\begin{corollary}
At the singular points of a surface of revolution all geodesics become tangent to
the direction of the axis projection $l_z$.
\end{corollary}

\proof
If the geodesic is ``vertical,'' i.e., $l_\gamma=l_z$, it stays
vertical forever, and the statement is evident. If its initial velocity
is non-vertical, $l_\gamma\not=l_z$, then the cross-ratio is initially finite.
Therefore, when the geodesic hits a singular point, the cross-ratio still 
has a finite limit which is obtained from the Clairaut invariant 
along the geodesic. On the other hand, when approaching the ``tropic,''
the null direction $l_{null}$ tends to the axis projection $l_z$. This forces the 
tangent element $l_\gamma$ to approach $l_z$ as well.
\proofend

Finally, consider space-like geodesics. Depending on the initial velocity, 
they can either hit the tropic or stay away from it. It is interesting to compare
the latter with geodesics in the Riemannian case. The Riemannian Clairaut
theorem implies that a ``non-vertical'' geodesic  does not enter the regions 
with a too narrow neck on the surface of revolution. Indeed, the invariance
of $r\sin\alpha$ implies, for  $\alpha\not=0$, that $r$ cannot be too 
small, since $|\sin\alpha|$ is bounded above. 
It turns out that, on a Lorentz surface 
of revolution, the phenomenon is exactly the opposite: such geodesics 
do not enter the regions where the neck is too wide:

\begin{corollary}\label{nowideneck}
For any space-like geodesic, there is an upper bound $K$ 
such that the geodesic stays in the region $|r|\le K$.
\end{corollary}

\proof
For space-like geodesics, $|{\rm cr}|< 1$. Then $|1-{\rm cr}^2|\le 1$ and
the conservation of $(1-{\rm cr}^2)/r^2$ implies that $|r^2|$ must be 
bounded above for any such geodesic.
\proofend

A surface of revolution is an example of a warped product.
A nice alternative description of the geodesics 
via Maupertuis' principle in the general context of warped products 
is given in \cite{Ze}. The dichotomy of the Riemannian and Lorentzian cases,
shown in Corollary \ref{nowideneck}, also follows from this description
of geodesics as particles moving in the potentials differed by sign.


\subsection{Analogs of the Jacobi and the Chasles theorems} \label{psellipt}
An ellipsoid with distinct axes in Euclidean space
\begin{equation} \label{psoid}
\frac{x_1^2}{a_1^2} + \frac{x_2^2}{a_2^2} +\dots + \frac{x_n^2}{a_n^2}=1
\end{equation}
gives rise to the confocal family of quadrics
$$
\frac{x_1^2}{a_1^2+\lambda} + \frac{x_2^2}{a_2^2+\lambda} +\dots +
\frac{x_n^2}{a_n^2+\lambda}=1.
$$

The Euclidean theory of confocal quadrics comprises the following theorems: through a generic point in space there pass $n$ confocal quadrics, and they are pairwise orthogonal at this point (Jacobi);\footnote{The respective values of $\lambda$ are called the elliptic coordinates of the point.} a generic line is tangent to $n-1$ confocal quadrics whose tangent hyperplanes at the points of tangency with the line are pairwise orthogonal (Chasles); and the tangent lines to a geodesic on an ellipsoid are tangent to fixed $n-2$ confocal quadrics (Jacobi-Chasles) -- see \cite{Mo1,Mo2,A-G}.

We shall construct a pseudo-Euclidean analog of this theory and adjust the proofs accordingly.

Consider pseudo-Euclidean space $V^n$ with signature $(k,l),\ k+l=n$ and let $E: V\to V^*$ be the self-adjoint operator such that $\langle x,x\rangle =E(x) \cdot x$ where dot denotes the pairing between vectors and covectors. Let $A: V\to V^*$ be a positive-definite self-adjoint operator defining an ellipsoid $A(x)\cdot x=1$. Since $A$ is positive-definite, both forms can be  simultaneously reduced to principle axes.\footnote{In general, it is not true that a pair of quadratic forms can be simultaneously reduced to principle axes. The simplest example in the plane is $x^2-y^2$ and $xy$.} 
We assume that $A={\rm diag}(a_1^2,\dots,a_n^2)$ and $E={\rm diag}(1,\dots,1,-1,\dots,-1)$. An analog of the confocal family is the following ``pseudo-confocal" family of quadrics  $Q_{\lambda}$
\begin{equation} \label{psconf}
\frac{x_1^2}{a_1^2+\lambda} + \frac{x_2^2}{a_2^2+\lambda} +\dots + \frac{x_k^2}{a_k^2+\lambda} + \frac{x_{k+1}^2}{a_{k+1}^2-\lambda}+\dots + \frac{x_n^2}{a_n^2-\lambda}=1
\end{equation}
where $\lambda$ is a real parameter or, in short, $(A+\lambda E)^{-1}(x)\cdot x=1$. 

The following result is a pseudo-Euclidean version of the Jacobi theorem.

\begin{theorem} \label{pelcoord}
Through every generic point $x\in V$ there pass either $n$ or $n-2$ quadrics from the pseudo-confocal family (\ref{psconf}). In the latter case, all quadrics have different topological types and in the former two of them have the same type. The quadrics are pairwise orthogonal at point $x$.
\end{theorem}

\begin{figure}[hbtp]
\centering
\includegraphics[width=3.5in]{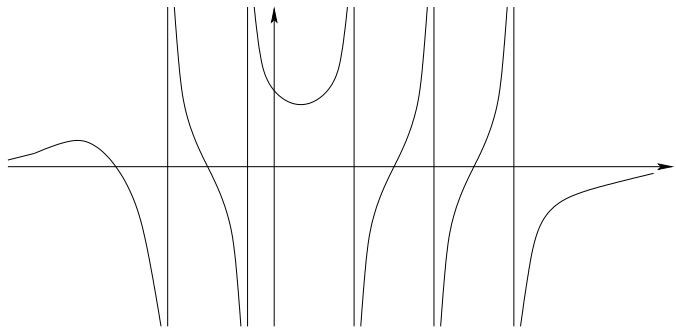}
\caption{The graph of the function $f(\lambda)$ for a time-like point $x$}
\label{graph}
\end{figure}

\proof  Given a point $x$, we want to find $\lambda$ satisfying  equation  (\ref{psconf}), which reduces to a polynomial in $\lambda$ of degree $n$.  Denote by $f(\lambda)$ the function on the left-hand-side of (\ref{psconf}). This function has
poles at $\lambda=-a_1^2,\dots,-a_k^2,a_{k+1}^2,\dots,a_n^2$.
At every negative pole $f(\lambda)$ changes sign from negative to positive, and at every positive pole from positive to negative. Let us analyze the behavior of $f(\lambda)$ as $\lambda \to \pm\infty$. One has:
$$
f(\lambda)= \frac{1}{\lambda} \langle x,x\rangle - \frac{1}{\lambda^2} \sum_{i=1}^n a_i^2x_i^2 +O\left(\frac{1}{\lambda^3} \right),
$$
hence if $x$ is not light-like then the sign of $f(\lambda)$ at $+\infty$ is equal, and at $-\infty$ opposite, to that of  $\langle x,x\rangle$, whereas if $x$ is  light-like then $f(\lambda)$ at $\pm\infty$ is negative. The graph of the function $f(\lambda)$ in the case $\langle x,x\rangle < 0$ is shown in  Figure \ref{graph}. Thus $f(\lambda)$ assumes value 1 at least $k-1$ times for negative $\lambda$ and at least $l-1$ times for positive ones. Being a polynomial of degree $n$, the number of roots is not greater than $n$. 

Note that the topological type of the quadric changes each time  that $\lambda$ passes through a pole of $f(\lambda)$. It follows that if there are $n-2$ quadrics passing through $x$ then they all have different topological types, and the  ellipsoid (corresponding to $\lambda=0$) is missing. On the other hand, if  there are $n$ quadrics passing through $x$ then two of them have the same  topological type and there are $n-1$ different  types altogether. Note, in particular, that if $x$ lies on the original ellipsoid then there are $n$ quadrics passing through it.

To prove that $Q_{\lambda}$ and $Q_{\mu}$ are orthogonal to each other at $x$, consider their normal vectors (half the gradients of the left-hand-sides  of (\ref{psconf}) with respect to the pseudo-Euclidean metric)
$$
N_{\lambda}= \left(
\frac{x_1}{a_1^2+\lambda},  \frac{x_2}{a_2^2+\lambda},\dots ,\frac{x_k}{a_k^2+\lambda},
-\frac{x_{k+1}}{a_{k+1}^2-\lambda},\dots,-\frac{x_n}{a_n^2-\lambda}\right),
$$
and likewise for $N_{\mu}$. Then 
\begin{equation} \label{dotns}
\langle N_{\lambda},N_{\mu}\rangle = 
\sum_{i=1}^k \frac{x_i^2}{(a_i^2+\lambda)(a_i^2+\mu)}-\sum_{i=k+1}^n \frac{x_i^2}{(a_i^2-\lambda)(a_i^2-\mu)}.
\end{equation}
The difference of the left-hand-sides  of equations  (\ref{psconf}), taken for $\lambda$ and $\mu$, is equal to the right-hand-side of (\ref{dotns}) times $(\mu -\lambda)$, whereas the right-hand-side is zero. Thus $\langle N_{\lambda},N_{\mu}\rangle =0$.
\proofend

\begin{figure}[hbtp]
\centering
\includegraphics[width=3in]{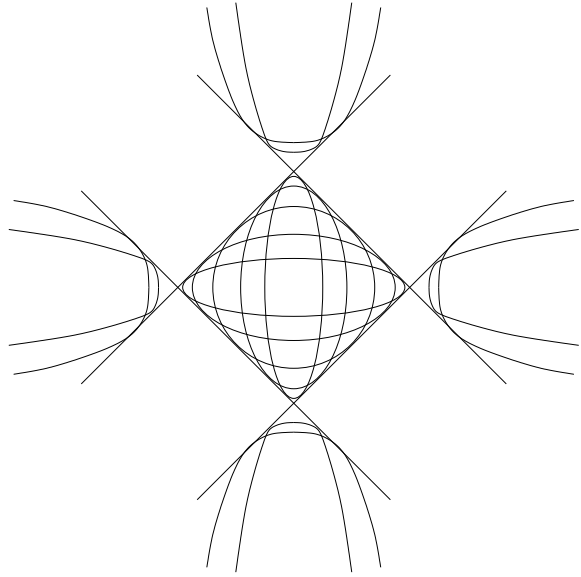}
\caption{A pseudo-confocal family of conics}
\label{domains}
\end{figure}

\begin{example} \label{exdimtwo}
{\rm Consider the simplest example, which we will study in detail 
in Section \ref{case} : $A={\rm diag}(1,1), E={\rm diag}(1,-1)$. Figure \ref{domains} depicts the partition of the Lorentz plane according to the number of conics from a pseudo-confocal family passing through a point: the boundary consists of the lines $|x\pm y|=\sqrt{2}$. }
\end{example}

\begin{problem} \label{bdry}
{\rm It is interesting to describe the topology of the partition of $V$ according to the number of quadrics from the family (\ref{psconf}) passing through a point. In particular, how many connected components are there?}
\end{problem}

Next, consider a pseudo-Euclidean version of the Chasles theorem.

\begin{theorem} \label{Chasles}
A generic space- or time-like line $\ell$  is tangent to either  $n-1$ or $n-3$, and a generic light-like line  to either  $n-2$ or $n-4$, quadrics from the family (\ref{psconf}). The tangent hyperplanes to these quadrics at the tangency points with $\ell$ are pairwise orthogonal. 
\end{theorem}

\proof Let $v$ be a vector spanning $\ell$. Suppose first that $v$ is space- or time-like. Project $V$ along $\ell$ on the orthogonal complement $U$ to $v$. A quadric determines a hypersurface in
this $(n-1)$-dimensional space, the set of critical values of its  projection (the
apparent contour). If one knows that these hypersurfaces also constitute a family (\ref{psconf})   of quadrics, the statement will follow from Theorem \ref{pelcoord}.

Let $Q\subset V$ be a smooth star-shaped hypersurface and let $W\subset V^*$ be the annihilator of $v$.  Suppose that a line parallel to $v$ is tangent to $Q$ at point $x$. Then the tangent hyperplane $T_x Q$ contains $v$. Hence the respective covector $y\in V^*$ from the polar dual hypersurface $Q^*$ lies in $W$. Thus polar duality takes the points of tangency of $Q$ with the lines parallel to $v$ to the intersection of the dual hypersurface $Q^*$ with the
hyperplane $W$.

On the other hand, $U=V/(v)$ and $W=(V/(v))^*$. Therefore the apparent contour of $Q$ in $U$ is polar dual to $Q^*\cap W$. If $Q$ belongs to the family of quadrics (\ref{psconf}) then $Q^*$ belongs to the pencil $(A+\lambda E) y\cdot y=1$. The intersection of a pencil 
with a hyperplane is a pencil of the same type (with the new $A$ positive definite and the new $E$ having signature $(k-1,l)$ or $(k,l-1)$, depending on whether $v$ is space- or time-like). It follows that the  polar dual family of quadrics, consisting of the apparent contours, is of the type (\ref{psconf}) again, as needed.

Note that, similarly to the proof of Theorem \ref{pelcoord}, if $\ell$ is tangent to the original ellipsoid then it is tangent to $n-1$ quadrics from the family (\ref{psconf}).

If $v$ is light-like then we argue similarly. We choose as the ``screen" $U=V/(v)$ any hyperplane transverse to $v$. The restriction of $E$ to $W$ is degenerate: it has 1-dimensional kernel and its signature is $(k-1,l-1,1)$. The family of quadrics, dual to the restriction of the pencil to $W$, is given by the formula 
$$
\frac{x_1^2}{b_1^2+\lambda} + \frac{x_2^2}{b_2^2+\lambda} +\dots + \frac{x_{k-1}^2}{b_{k-1}^2+\lambda} + \frac{x_{k+1}^2}{b_{k+1}^2-\lambda}+\dots + \frac{x_{k+l-1}^2}{b_{k+l-1}^2-\lambda}
=1-\frac{x_{k}^2}{b_{k}^2}
$$
which is now covered by the $(n-1)$-dimensional case  of Theorem \ref{pelcoord}.
\proofend

Note that in Example \ref{exdimtwo} a generic light-like line is tangent to no conic, whereas the four exceptional light-like lines $|x\pm y|=\sqrt{2}$ are tangent to infinitely many ones.


\subsection{Complete integrability} \label{inv}

The following theorem is a pseudo-Euclidean analog of the Jacobi-Chasles theorem.\footnote{See also \cite{GKT}, a follow-up to the present paper, devoted to the case study of the geodesic flow on the ellipsoid in 3-dimensional pseudo-Euclidean space.} 

\begin{theorem} \label{JCh}
1) The tangent lines to a fixed space- or time-like (respectively, light-like) geodesic on a quadric in  pseudo-Euclidean space $V^n$ are tangent to $n-2$ (respectively, $n-3$) other fixed quadrics from the pseudo-confocal family  (\ref{psconf}).\\
2) A space- or time-like (respectively, light-like) billiard trajectory in a quadric in pseudo-Euclidean space  $V^n$ remains tangent to $n-1$ (respectively, $n-2$) fixed quadrics from the family  (\ref{psconf}).\\
3) The sets  of space- or time-like  oriented lines in pseudo-Euclidean  space $V^n$, tangent to  $n-1$  fixed  quadrics from the family  (\ref{psconf}), are  Lagrangian  submanifolds in the spaces ${\cal L}_{\pm}$. The set of light-like oriented lines, tangent to  $n-2$  fixed  quadrics from the family  (\ref{psconf}), is a codimension $n-2$  submanifold in   ${\cal L}_0$ foliated by codimension one Legendrian submanifolds. 
\end{theorem}

\proof Let $\ell$ be a tangent line at point $x$ to a geodesic on the quadric $Q_0$ from the pseudo-confocal family  (\ref{psconf}). By Theorem \ref{Chasles}, $\ell$ is tangent to $n-2$ (or $n-3$, in the light-like case) quadrics from this family. Denote these quadrics by $Q_{\lambda_j},\ j=1,\dots, n-2$.

Let $N$ be a normal vector to $Q_0$ at point $x$. Consider an infinitesimal 
rotation of the tangent line $\ell$ along the geodesic. Modulo infinitesimals of the second order, this line rotates in the 2-plane generated by $\ell$ and $N$. By Theorem \ref{Chasles}, the tangent hyperplane to each $Q_{\lambda_j}$ at its tangency point with $\ell$ contains the  vector $N$. Hence, modulo infinitesimals of the  second order, the line $\ell$ remains tangent to every $Q_{\lambda_j}$, and therefore remains tangent to each one of them.

The billiard flow inside an ellipsoid in $n$-dimensional space is the  limit case of the geodesic flow on an ellipsoid in $(n+1)$-dimensional space, whose  minor axis goes  to zero. Thus the second statement follows from the first one.

Now we prove the third statement. Consider first the case of space- or time-like 
lines ${\cal L}_{\pm}$. Let $\ell$ be a generic oriented line tangent to quadrics $Q_{\lambda_j},\ j=1,\dots,n-1$, from the family  (\ref{psconf}). Choose smooth functions $f_j$ defined in neighborhoods of the tangency points of $\ell$ with $Q_{\lambda_j}$ 
in $V^n$ whose level hypersurfaces are the quadrics from the family  (\ref{psconf}). Any line  $\ell'$ close to $\ell$ is tangent to a close quadric $Q_{\lambda'_j}$. Define the function $F_j$ on the space of oriented lines whose value at $\ell'$ is the (constant) value of $f_j$ on $Q_{\lambda'_j}$.

We want to show that $\{F_j,F_k\}=0$ where the Poisson bracket is taken with respect to the symplectic structure defined in Section \ref{str}. Consider the value $dF_k({\rm sgrad}\ F_j)$ at $\ell$. The vector field ${\rm sgrad}\  F_j$ is tangent  to the characteristics of the hypersurface $F_j$ = const, that is,  the hypersurface  consisting  of the lines, tangent to $Q_{\lambda_j}$. According to Theorem \ref{geoline}, these characteristics consist of the lines, tangent  to a fixed geodesic  on $Q_{\lambda_j}$. According to statement 1 of the present theorem, these lines are tangent to $Q_{\lambda_k}$, hence $F_k$ does not change along the flow of ${\rm sgrad}\  F_j$. Thus $dF_k({\rm sgrad}\ F_j)=0$, as claimed. 

Finally, in the  light-like case, consider the homogeneous symplectic manifold $P^{2n-2}$ of scaled light-like lines whose quotient is ${\cal L}_0$. Then, as before, we have homogeneous of degree zero, Poisson-commuting functions $F_j,\ j=1,\dots,n-2$, on $P$. Therefore a level submanifold $M^n=\{F_1=c_1,\dots,F_{n-2}=c_{n-2}\}$ is coisotropic: the symplectic orthogonal complement to $TM$ in $TP$ is contained in $TM$. 
The commuting vector fields sgrad $F_j$ define an action of the Abelian group $\R^{n-2}$ on $M$ whose orbits are isotropic submanifolds.
Furthermore, $M$ is invariant under  the Euler vector field $E$ that preserves the foliation  on isotropic submanifolds. Hence the quotient by $E$ is a codimension $n-2$  submanifold in ${\cal L}_0$ foliated by codimension one Legendrian submanifolds.
\proofend

\begin{example} \label{3ell}
{\rm Let $\gamma$ be a geodesic on a generic ellipsoid $Q_0$ in 3-dimensional Lorentz space and let $x$ be a point of $\gamma$. Then, upon each return to point $x$, the curve $\gamma$ has one of  at most two possible directions (a well known property in the Euclidean case). 

Indeed, if $\gamma$ is not light-like then the tangent lines to $\gamma$ are tangent to a fixed pseudo-confocal quadratic surface, say $Q_1$. The intersection of the tangent plane $T_x Q_0$ with $Q_1$ is a conic, and there are at most two tangent lines from $x$ to this conic. If $\gamma$ is light-like then its direction at point $x$ is in the kernel of the restriction of the metric to $T_x Q_0$ which consists of at most two lines. 
}
\end{example}

\begin{remark} \label{comm}
{\rm The functions $F_1,\dots, F_{n-1}$ from the proof of Theorem \ref{JCh} can be considered as functions on the tangent bundle $TM$. The proof of  Theorem \ref{JCh} implies that these functions and the energy function $F_0(x,v)=\langle v,v\rangle$ pairwise Poisson commute with respect to the canonical symplectic structure on $T^*M$ (as usual, identified with $TM$ via the metric).

We now give explicit formulas for the integrals.; these formulas are modifications of the ones given in \cite{Mo1,Mo2}.

Write the pseudo-Euclidean metric as 
$$
\sum_{i=1}^n \tau_i dx_i^2\quad{\rm with}\quad \tau_1=\dots=\tau_k=1,\ \tau_{k+1}=\dots=\tau_n=-1.
$$
 Let $v=(v_1,\dots,v_n)$ denote tangent vectors to the ellipsoid. 
Then the integrals are given by the formulas
$$
F_k=\frac{v_k^2}{\tau_k} +\sum_{i\neq k} \frac{(x_iv_k-x_kv_i)^2}{\tau_i a_k^2-\tau_k a_i^2},\ \ k=1,\dots,n.
$$
These integrals satisfy the relation $\sum F_k = \langle v,v\rangle$. 
One also has a modification of the Joachimsthal  integral (functionally dependent on the previous ones):
$$
J=\left( \sum_i \frac{x_i^2}{\tau_i a_i^4}\right) \left(\sum_j \frac{v_j^2}{a_j^2}\right).
$$
}
\end{remark}

\begin{remark} \label{alter}
{\rm 
Another approach to complete integrability of the billiard in the ellipsoid and 
the geodesic flow  on the ellipsoid in Euclidean space  is described in \cite{Ta6,Ta3}. In a nutshell, in the case of billiards, one constructs another symplectic form on the space of oriented lines invariant under the billiard map, and for the geodesic flow one constructs another metric on the ellipsoid, projectively equivalent to the Euclidean one: this means that their non-parameterized geodesics coincide. For geodesic flows, this  integrability mechanism was independently and simultaneously discovered by  Matveev and Topalov \cite{M-T1}.

In the present situation, this approach leads to the following result. We do not dwell on details. 

As before,  $Q$ is an ellipsoid in pseudo-Euclidean space $V^n$ given by the equation $A(x)\cdot x=1$ and  the scalar product is $\langle u,v\rangle=E(u)\cdot v$ where $A$ and $E$ are self-adjoint operators $V\to V^*$. We assume that $A={\rm diag}(a_1^2,\dots,a_n^2)$ and $E={\rm diag}(1,\dots,1,-1,\dots,-1)$, and denote the billiard map in $Q$ by $T$. Consider the interior of $Q$ as the projective, or Cayley-Klein, model of hyperbolic geometry and let $\Omega$ be the respective symplectic structure on the space of oriented lines (obtained by the standard symplectic reduction).

\begin{theorem} \label{Lob}
1) The symplectic structure $\Omega$ is invariant under $T$. \\
2) The restrictions of the metrics 
\begin{equation} \label{2met}
\langle dx,dx\rangle\quad {\rm and}\quad 
\frac{A(dx)\cdot dx}{\langle A(x),A(x)\rangle}
\end{equation}
on the ellipsoid $Q$ are projectively equivalent.
\end{theorem}
}
\end{remark}


\section{The geometry of the circle
billiard in the Lorentz plane} \label{case}

Below we show how the theorems of Section \ref{quad} 
work for a circle billiard in dimension 2. 
Although its integrability follows from the results above,
we made this section self-contained to emphasize 
the simplicity of the corresponding formulae.

Consider the plane with the metric $ds^2=dx dy$. Then a vector $(a,b)$ is orthogonal to $(a,-b)$. Let $D(a,b)=(b,a)$ be the linear operator identifying vectors and covectors via the metric.

Consider the circle $x^2+y^2=1$ and the billiard system inside it. There are four singular points: $(\pm1,0),(0,\pm1)$. The phase space consists of the oriented lines intersecting the circle and such that the impact point is not singular. The billiard map on light-like lines is 4-periodic. One also has two 2-periodic orbits, the diameters having slopes $\pm 1$.

Let $t$ be the cyclic coordinate on the circle. Let us characterize a line by the coordinates of its first and second intersection points with the circle, $(t_1, t_2)$. The billiard map  $T$ sends $(t_1,t_2)$ to $(t_2,t_3)$. 

\begin{theorem} \label{circlebil}
1) The map $T$ is given by the equation
\begin{equation} \label{forcirc}
\cot \left(\frac{t_2-t_1}{2}\right)  + \cot \left(\frac{t_2-t_3}{2}\right) =2\cot 2t_2.
\end{equation}
2) The area form is given by the formula
\begin{equation} \label{arcirc}
\omega=\frac{\sin \left((t_2-t_1)/2\right)}{|\sin (t_1+t_2)|^{3/2}}\ dt_1\wedge dt_2.
\end{equation}
3) The map is integrable: it has an invariant function
\begin{equation} \label{intcirc}
I=\frac{\sin \left((t_2-t_1)/2\right)}{|\sin (t_1+t_2)|^{1/2}}.
\end{equation}
4)  The lines containing the billiard segments, corresponding to a fixed value
$\lambda$ of the (squared) integral  $I^2$, are tangent to the  conic
\begin{equation} \label{conics}
x^2+y^2+2\lambda xy=1-\lambda^2,\ \ \lambda\in \R.
\end{equation}
These conics for different $\lambda$ are all tangent to the four lines -- two horizontal and two vertical --  tangent to the unit circle.
\end{theorem}

In principal axes (rotated $45^{\circ}$), the family (\ref{conics}) writes as
$$
\frac{x^2}{1-\lambda}+\frac{y^2}{1+\lambda}=1.
$$

Before going into the proof of Theorem \ref{circlebil} let us make some comments and illustrate the theorem by figures.

In the familiar case of the billiard inside an ellipse in the Euclidean plane, the billiard trajectories are tangent to the family of conics, confocal with the given ellipse, see, e.g., \cite{Ta5,Ta2}. These conics are either the confocal ellipses inside the elliptic billiard table or the confocal hyperbolas. The billiard map, restricted to an invariant curve of the integral, is described as follows: take a point $A$ on the boundary of an elliptic billiard table, draw a tangent line to the fixed confocal ellipse (or hyperbola -- for other values of the integral) until the intersection with the boundary ellipse at point $A_1$; take $A_1$ as the next point of the billiard orbit, etc.

\begin{figure}[hbtp]
\centering
\includegraphics[width=1.5in]{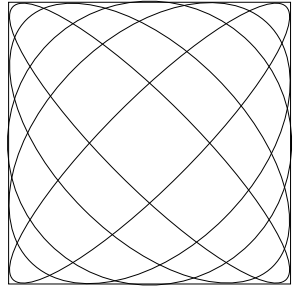}
\caption{Ellipses inscribed into a square}
\label{square}
\end{figure}

\begin{figure}[hbtp]
\centering
\includegraphics[width=1.5in]{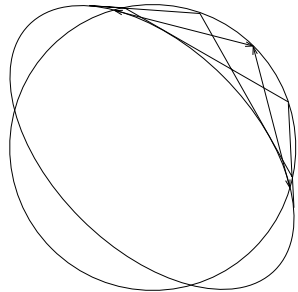}
\caption{Two 4-periodic trajectories on the invariant curve with the rotation number $1/4$: one orbit is a self-intersecting quadrilateral, and the other one consists of two segments, traversed back and forth}
\label{per4}
\end{figure}

In our case, for a fixed billiard trajectory inside the unit circle, 
there is a quadric inscribed into a $2\times 2$ square 
to which it is tangent, see the family of ellipses in figure \ref{square}. 
For instance, two 4-periodic trajectories in figure \ref{per4} are tangent to one and the same inscribed ellipse.

\begin{figure}[hbtp]
\centering
\includegraphics[width=2in]{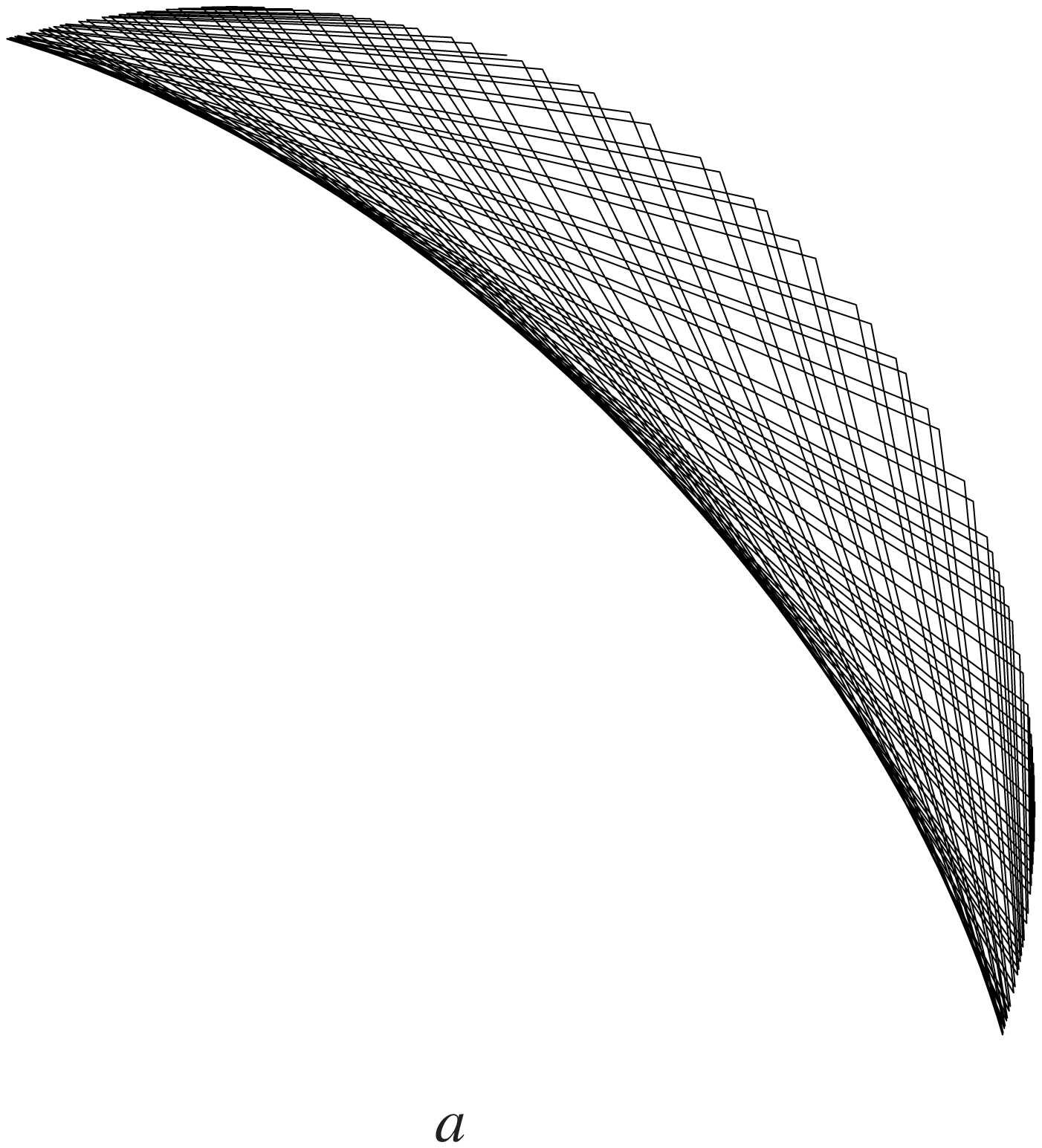}
\qquad
\includegraphics[width=2in]{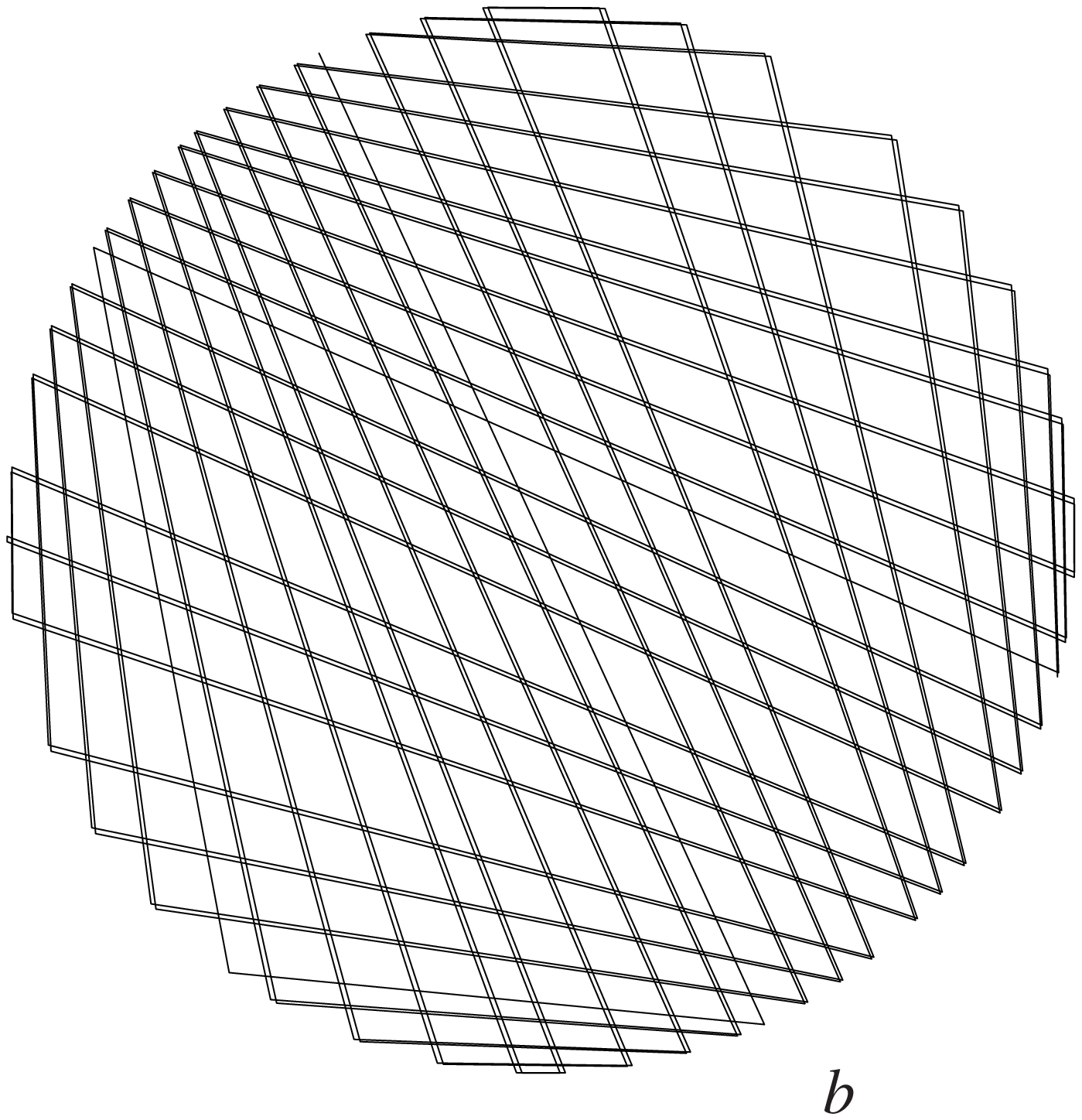}
\caption{Two orbits consisting of 100 billiard segments: the lines containing the billiard segments on the left are tangent to an ellipse, and those on the right to a hyperbola}
\label{100lines}
\end{figure}

Figures \ref{100lines}$a$ and \ref{100lines}$b$ 
depict two billiard orbits in the configuration 
space consisting of 100 time-like billiard segments. 
It is easy to recognize the inscribed ellipse as an envelope of the segments 
of the  billiard orbit in figure \ref{100lines}$a$: all the reflections occur on one of the two arcs of the circle outside of the ellipse, and the tangent line to the ellipse at its intersection point with the circle is the Lorentz normal to the circle at this point.

\begin{figure}[hbtp]
\centering
\includegraphics[width=2.5in]{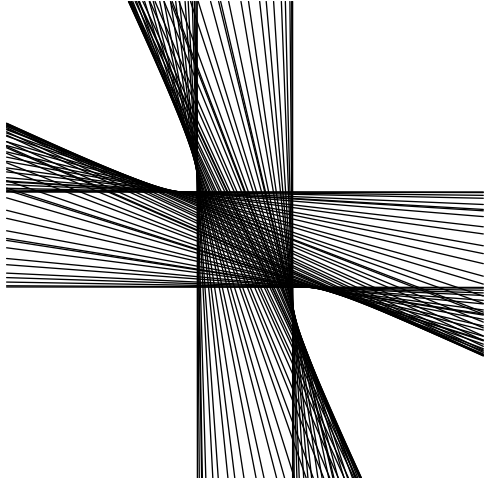}
\caption{Extensions of billiard chords tangent to a hyperbola}
\label{envhyp}
\end{figure}

The corresponding envelope is less evident for the orbit in figure \ref{100lines}$b$: 
extensions of the billiard chords have a hyperbola as their envelope, see
figure \ref{envhyp}.

\begin{figure}[hbtp]
\centering
\includegraphics[width=3in]{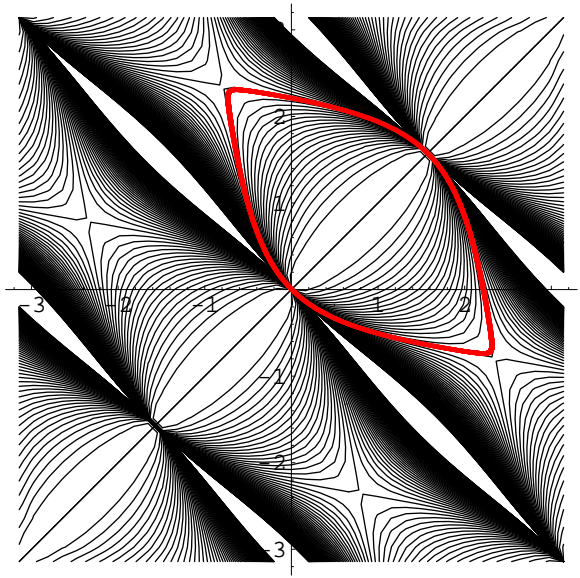}
\caption{Level curves of the integral $I$ in the $(t_1,t_2)$ coordinates} 
\label{folio}
\end{figure}

The level curves of the integral $I$ are shown in figure \ref{folio} depicting a $[-\pi,\pi]\times [-\pi,\pi]$ torus with coordinates $(t_1,t_2)$. The four hyperbolic singularities of the foliation $I$=const at points 
$$(3\pi/4,-\pi/4), (\pi/4,-3\pi/4),(-\pi/4,3\pi/4),(-3\pi/4,\pi/4)$$
correspond to two 2-periodic orbits of the billiard map; these orbits are hyperbolically unstable\footnote{As indicated by the hyperbolic crosses made by the level curves at these points.} (unlike the case of an ellipse in the Euclidean plane where the minor axis is a stable 2-periodic orbit). The white spindle-like regions surround the lines $t_1+t_2=\pi n,\ n\in \Z$; these lines correspond to the light-like rays. The four points 
$(0,0), (\pi/2,\pi/2), (-\pi/2,-\pi/2), (\pi,\pi)$ are singular: every level curve of the integral $I$ pass through them. 

\begin{figure}[hbtp]
\centering
\includegraphics[width=1.5in]{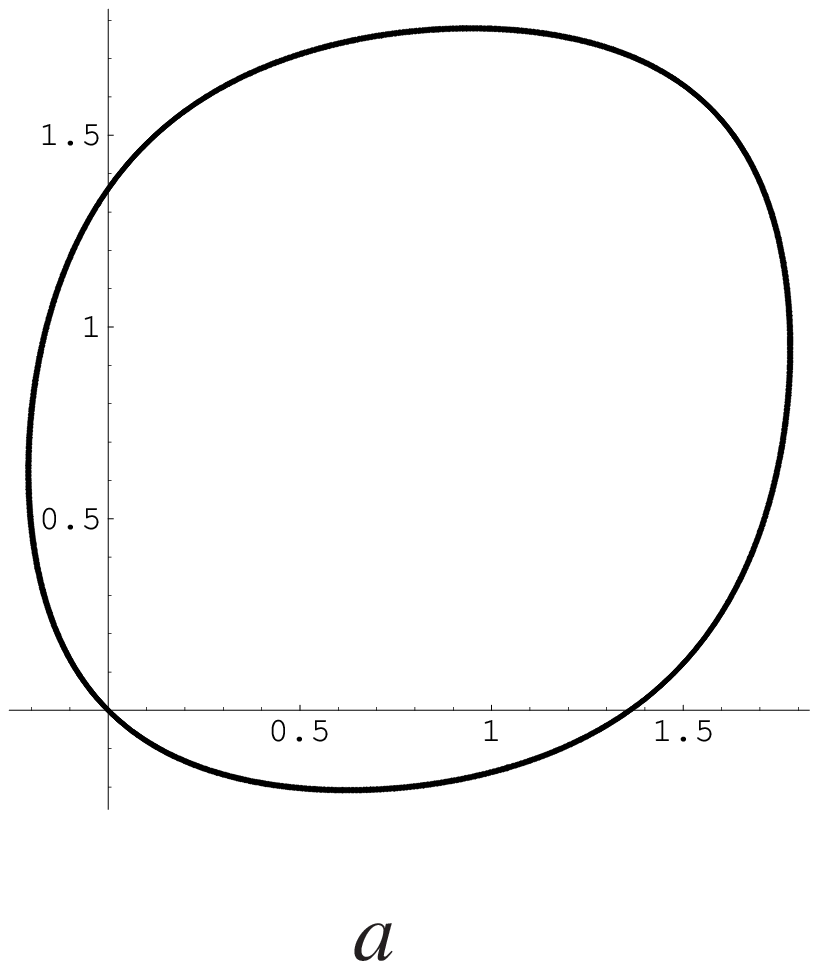}
\qquad
\includegraphics[width=1.5in]{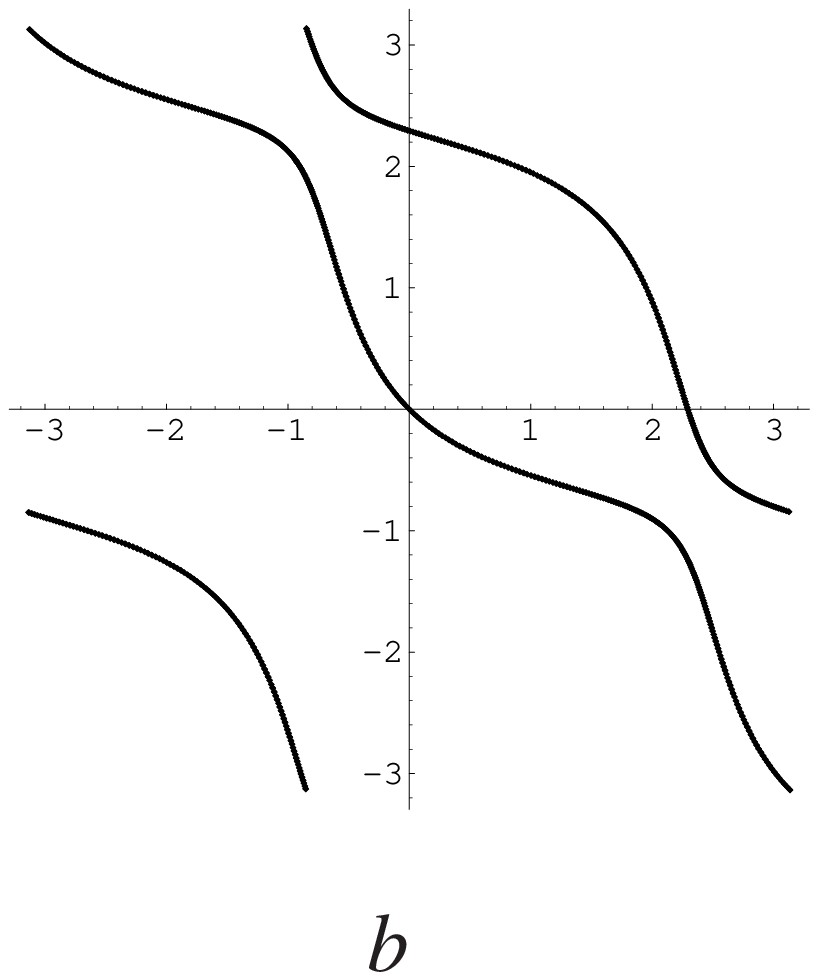}
\caption{Two invariant curves (the right one consists of two components, the phase points ``jump" from one component to the other)}
\label{invc}
\end{figure}

Figures \ref{invc}$a$ and \ref{invc}$b$ show 
two topologically different invariant curves, and 
figures \ref{100lines}$a$ and \ref{100lines}$b$  discussed above depict 
two billiard orbits  in the configuration space corresponding respectively to those two invariant curves.

Note also that it follows from the Poncelet porism (see, e.g., \cite{B-K-O-R}) that  if some point of an invariant circle is periodic then all points of this   invariant circle are periodic with the same period, cf. figure \ref{per4}.

Now we shall prove Theorem \ref{circlebil}.

\begin{figure}[hbtp]
\centering
\includegraphics[width=2in]{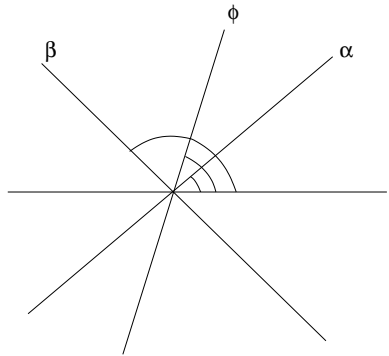}
\caption{Harmonic quadruple of lines given by angles}
\label{harmo}
\end{figure}

\proof
We use another criterion for harmonicity of a quadruple of lines, similar to (\ref{harmcr}). Consider four concurrent lines, and let $\alpha,\phi,\beta$ be the angles made by three of them with the fourth, see figure \ref{harmo}. Then the lines are harmonic if and only if
\begin{equation} \label{harm}
\cot \alpha + \cot \beta = 2 \cot \phi.
\end{equation}

In our situation, the billiard curve is $\gamma(t)=(\cos t, \sin t)$. The tangent vector is $\gamma'(t)=(-\sin t,\cos t)$ and the normal is $(\sin t, \cos t)$. Consider the impact point $t_2$. By elementary geometry, the rays $(t_2,t_1)$ and $(t_2,t_3)$ make the angles $(t_1-t_2)/2$ and $(t_3-t_2)/2$ with the tangent line at $\gamma(t_2)$, and the normal makes the angle $\pi-2t_2$ with this tangent line. Then equation (\ref{harm}) becomes (\ref{forcirc}).

It is straightforward to compute the area form from Lemma \ref{symstr2} in the $(t_1,t_2)$-coordinates; the result (up to a constant factor) is (\ref{arcirc}). 

We shall give two proofs that $I$ is an integral. First, our Lorentz billiard is a particular case of a projective billiard in a circle. It is proved in \cite{Ta1} that every such billiard map has an invariant area form 
\begin{equation} \label{invform}
\Omega= \frac{1}{\sin^2 ((t_2-t_1)/2)}\ dt_1\wedge dt_2.
\end{equation}
(This form is the symplectic structure on the space of oriented lines for the projective -- or Klein-Beltrami -- model of hyperbolic geometry inside the unit disc, see Remark \ref{alter}.) Thus $T$ has two invariant area forms, and (the cube root of) their ratio is an invariant function.

The second proof imitates a proof that the billiard inside an ellipse in the Euclidean plane is integrable, see \cite{Ta2}. Let us restrict attention to space-like lines. Assign to a line its first intersection point with the circle, $q$, and the unit vector along the line, $v$. Then $\langle D(q),q\rangle=1$ and $\langle v,v\rangle=1$. We claim that $I=\langle D(q),v\rangle$ is invariant under the billiard map. 

The billiard map is the composition of the involutions $\tau$ and $\sigma$, 
see proof of Theorem \ref{invstr}. It turns out that each involution 
changes the sign of $I$. 

Indeed, $\langle D(q)+D(q_1),q_1-q\rangle=0$ since $D$ is self-adjoint. Since $v$ is collinear with $q_1-q$, we have: $\langle D(q)+D(q_1),v\rangle=0$, and hence $I$ is odd with respect to  $\tau$. 

Since the circle is given by the equation $\langle D(q),q\rangle=1$, the normal at point $q_1$ is $D(q_1)$. By definition of the billiard reflection, the vector $v+v_1$ is collinear with the normal at $q_1$, hence  $\langle D(q_1),v\rangle = -\langle D(q_1),v_1\rangle$. Thus $I$ is odd with respect to $\sigma$ as well. The invariance of $I=\langle D(q),v\rangle$ follows, and it is straightforward to check that, in the $(t_1,t_2)$-coordinates, this integral equals (\ref{intcirc}). 

\begin{figure}[hbtp]
\centering
\includegraphics[width=4in]{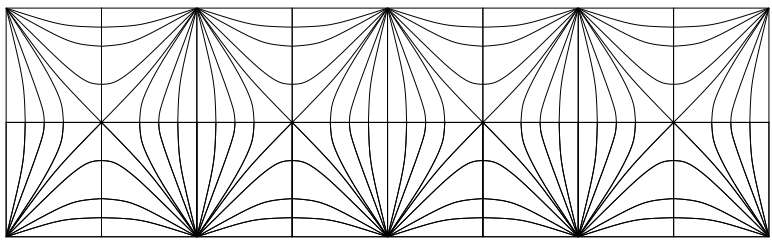}
\caption{Level curves of the integral $I$ in the $(\alpha,p)$-coordinates} 
\label{phspace}
\end{figure}

To find the equation of the envelopes and prove {\it 4)} we first rewrite 
the integral $I$ in the standard, Euclidean, coordinates $(p,\alpha)$
in the space of lines: $p$ is the signed length of the perpendicular 
from the origin to the line and $\alpha$ the direction of this perpendicular,
see \cite{San}. One has
\begin{equation} \label{pal}
\alpha=\frac{t_1+t_2}{2},\ p=\cos \left(\frac{t_2-t_1}{2}\right).
\end{equation}
Fix a value of the integral $I$ by setting
$$
\frac{\sin^2 ((t_2-t_1)/2)}{\sin(t_1+t_2)} =\lambda.
$$
It follows from (\ref{pal}) that $1-p^2=\lambda \sin 2\alpha$, and hence $p=\sqrt{1-\lambda \sin 2\alpha}$. (See figure  \ref{phspace} which shows the level curves of the (squared) integral $I^2=(1-p^2)/\sin 2\alpha$ in the $(\alpha,p)$-coordinates.) We use $\alpha$ as a coordinate on the level curve corresponding to a fixed value of $\lambda$, and $p$ as a function of $\alpha$ (this function depends on $\lambda$ as a parameter).

The envelope of a 1-parameter family of lines given by a function $p(\alpha)$ is the curve
$$
(x(\alpha),y(\alpha))=p(\alpha) (\cos \alpha,\sin \alpha)+p'(\alpha)(-\sin \alpha,\cos \alpha),
$$
see \cite{San}. 
In our case, we obtain the curve
$$
(x(\alpha),y(\alpha))=(1-\lambda \sin 2\alpha)^{-1/2} (\cos \alpha - \lambda \sin \alpha,\sin \alpha - \lambda \cos \alpha).
$$
It is straightforward to check that this curve satisfies equation (\ref{conics}).

It is also clear that the conics (\ref{conics}) are tangent to the lines $x=\pm1$ and $y=\pm 1$. Indeed, if, for example, $y= 1$ then the left hand side of (\ref{conics}) becomes $(x+ \lambda)^2+1-\lambda^2$, and equation (\ref{conics}) has a multiple root $x_{1,2}=-\lambda$.
\proofend

\begin{remark} \label{altpf}
{\rm Yet another proof of the integrability of the Lorentz billiard inside a circle can be deduced from the duality (the skew hodograph transformation) between Minkowski billiards discovered in \cite{GT}. This duality trades the shape of the billiard table for that of the unit (co)sphere of the metric. In our case, the billiard curve is a circle and the unit sphere of the metric is a hyperbola; the dual system is the usual, Euclidean billiard ``inside'' a hyperbola.
}
\end{remark}

\bigskip


\end{document}